\documentclass[preprint,11pt]{elsarticle}
\usepackage{amssymb}
\usepackage{amsmath}
\usepackage{amsfonts}
\usepackage{amssymb}
\usepackage{indentfirst,latexsym,bm}

\usepackage{tikz}
\usetikzlibrary{calc}

\usepackage{amsthm}
\usepackage{color,xcolor}
\usepackage[all]{xy}
\input{mathrsfs.sty}
\newtheorem{theorem}{Theorem}[section]
\newtheorem{lemma}[theorem]{Lemma}
\newtheorem{corollary}[theorem]{Corollary}
\newtheorem{proposition}[theorem]{Proposition}

\theoremstyle{definition}
\newtheorem{definition}{Definition}[section]

\theoremstyle{definition}
\newtheorem{example}{Example}[section]

\theoremstyle{remark}
\newtheorem{remark}{Remark}[section]
\theoremstyle{question}

\theoremstyle{problem}
\newtheorem{problem}{Problem}[section]

\numberwithin{equation}{section}
\def\be{\begin{equation}}
\def\ee{\end{equation}}

\journal{XXX}

\begin{document}

\begin{frontmatter}

%% Title, authors and addresses

%% use the tnoteref command within \title for footnotes;
%% use the tnotetext command for the associated footnote;
%% use the fnref command within \author or \address for footnotes;
%% use the fntext command for the associated footnote;
%% use the corref command within \author for corresponding author footnotes;
%% use the cortext command for the associated footnote;
%% use the ead command for the email address,
%% and the form \ead[url] for the home page:
%%
%% \title{Title\tnoteref{label1}}
%% \tnotetext[label1]{}
%% \author{Name\corref{cor1}\fnref{label2}}
%% \ead{email address}
%% \ead[url]{home page}
%% \fntext[label2]{}
%% \cortext[cor1]{}
%% \address{Address\fnref{label3}}
%% \fntext[label3]{}

\title{Some applications of the matched projections of idempotents}
\author[shnu]{Xiaofeng Zhang}
\ead{xfzhang8103@163.com}
\author[shnu]{Xiaoyi Tian}
\ead{tianxytian@163.com}
\author[shnu]{Qingxiang Xu}
\ead{qingxiang\_xu@126.com}
%\fntext[fn1]{Partially supported by the
%National Natural Science Foundation of China (11971136).}
\address[shnu]{Department of Mathematics, Shanghai Normal University, Shanghai 200234, PR China}

\begin{abstract}For every idempotent $Q$ on a Hilbert space $H$, the matched projection $m(Q)$ is a well-established concept. This paper explores several applications of the matched projections. The first application addresses the distances from projections on $H$ to a given idempotent $Q$. Using $m(Q)$, a complete characterization of these distances is established, covering the minimum, maximum, and intermediate values. The second application focuses on the $C^*$-algebra $C^*\{Q\}$ generated by a single non-projection idempotent $Q$. A new $4\times 4$ block matrix representation of $Q$, induced by $m(Q)$, yields novel formulas for $Q$, leading to a full characterization of all elements in $C^*\{Q\}$ via explicit $4\times 4$ block matrices. Furthermore, for each $r>1$, a family of universal $r$-idempotents is introduced. These idempotents possess a universal property distinct from known properties of projection pairs. Some necessary and sufficient conditions are provided for such universal $r$-idempotents. The third application presents new characterizations of the numerical ranges. An operator  version of the elliptical range theorem is established. Using a general non-projection idempotent $Q$ and its matched projection $m(Q)$, a non-quadratic operator is constructed, and its numerical range is described in detail. Additionally, another operator is introduced whose  numerical range closure is not an elliptical disk, and the numerical range itself is neither closed nor open.
\end{abstract}

\begin{keyword}Idempotent, projection, matched projection, operator distance, numerical range
\MSC 47A05, 46L05, 47A12, 15A60

%% MSC codes here, in the form: \MSC code \sep code

%% or \MSC[2008] code \sep code (2000 is the default)

\end{keyword}

\end{frontmatter}

%%
%% Start line numbering here if you want
%%
% \linenumbers

%% main text
%\tableofcontents
\tableofcontents
\section{Introduction}\label{sec:Intro}
Throughout this paper, $H$ and $K$ denote non-zero Hilbert spaces, and $\mathbb{B}(H,K)$ denotes the set of all bounded linear operators from $H$ to $K$.
If $H=K$, we write $\mathbb{B}(H)$ instead of $\mathbb{B}(H,H)$. For any $T\in \mathbb{B}(H)$, the symbol $\sigma_{\mathfrak{A}}(T)$ denotes the spectrum of $T$ with respect to a $C^*$-subalgebra $\mathfrak{A}$ of  $\mathbb{B}(H)$, and we abbreviate $\sigma_{\mathbb{B}(H)}(T)$ as $\sigma(T)$. Let $\sigma_p(T)$ denote the point spectrum of $T$, i.e., the set of its eigenvalues. The identity operator on $H$ is denoted by $I_H$, or simply $I$ when no confusion arises.

 For any $T\in \mathbb{B}(H,K)$, its  range, null space and adjoint operator are denoted by $\mathcal{R}(T)$, $\mathcal{N}(T)$ and $T^*$, respectively. Let $|T|$ denote the square root of $T^*T$. An operator  $Q\in \mathbb{B}(H)$ is called an idempotent if $Q^2=Q$; if in addition $Q^*=Q$, then $Q$ is  a projection.
 An idempotent that is not a projection is called a non-projection idempotent, while a projection $P$ is termed non-trivial if  $P\ne I$ and $P\ne 0$.
 Note that in some literature, such idempotents are also called skew projections \cite{BS,Ovchinnikov}.

For any $T\in \mathbb{B}(H,K)$ and any closed subspace $M$ of $H$, denote by  $T|_M$  the restriction of $T$ to $M$, and by
$P_M$ the projection from $H$ onto $M$.
Let $\mathbb{N}$ and $\mathbb{C}$ denote the sets of positive integers and complex numbers, respectively. Standard matrix notation is adopted: in particular, $M_n(\mathbb{C})$ denotes the set of $n\times n$ complex matrices,  and $I_n$ and $0_n$ represent  the identity and zero matrices  in $M_n(\mathbb{C})$, respectively.
For any topological space $\Omega$, let $C(\Omega)$ denote the space of complex-valued continuous functions on $\Omega$.

For Hilbert spaces $H_1$ and $H_2$, their sum $H_1\oplus
H_2$ is the Hilbert space defined by $$H_1\oplus
H_2=\left\{\binom{h_1}{h_2}:h_1\in H_1, h_2\in H_2\right\},$$ equipped with the inner-product
$$\left\langle\binom{h_1}{h_2},\binom{h_1'}{h_2'}\right\rangle=\langle h_1,h_1' \rangle + \langle h_2,h_2' \rangle,\quad \text{for all $h_i,h_i' \in H_i,i=1,2$}.$$

Given a non-trivial projection $P\in\mathbb{B}(H)$ (i.e.,\,$P\ne I$ and $P\ne 0$), define the unitary operator $U_P: H\to \mathcal{R}(P)\oplus \mathcal{N}(P)$ by
\begin{equation}\label{equ:unitary operator induced by P}U_Ph=\binom{Ph}{(I-P)h}, \quad h\in H.
\end{equation}
Its inverse $U_P^*$ acts as
\begin{equation*}\label{expression of inverse of U P}U_P^{*}\binom{h_1}{h_2}=h_1+h_2,\quad h_1\in\mathcal{R}(P), h_2\in\mathcal{N}(P).
\end{equation*}
For any  $T\in \mathbb{B}(H)$,  we have  (see e.g.,\,\cite[Section~2.2]{Xu-Wei-Gu})
 \begin{equation*}
U_PTU_P^{*}=\left(
            \begin{array}{cc}
              PTP|_{\mathcal{R}(P)} & PT(I-P)|_{\mathcal{N}(P)} \\
              (I-P)TP|_{\mathcal{R}(P)} &(I-P)T(I-P)|_{\mathcal{N}(P)}\\
            \end{array}
          \right).
\end{equation*}
In particular,
 \begin{equation*}\label{equ:block matrix P}
U_PPU_P^{*}=\left(
            \begin{array}{cc}
              I_{\mathcal{R}(P)} & 0 \\
              0 &0\\
            \end{array}
          \right),\quad U_P(I-P)U_P^{*}=\left(
            \begin{array}{cc}
              0 & 0 \\
              0 &I_{\mathcal{N}(P)}\\
            \end{array}
          \right).
\end{equation*}
Conversely, any operator $X\in \mathbb{B}\big(\mathcal{R}(P)\oplus \mathcal{N}(P)\big)$ with the block form
$$X=\left(\begin{array}{ccc} X_{11}&
X_{12}\\X_{21}&X_{22}\end{array}\right)$$
is mapped back to $\mathbb{B}(H)$ via
\begin{equation*}U_P^{*}XU_P=(X_{11}+X_{21})P+(X_{12}+X_{22})(I-P).\end{equation*}

Idempotents and projections are two fundamental and closely related concepts in operator theory and linear algebra, which have been the subject of extensive research due to their theoretical significance and practical applications. Despite the substantial body of existing work, certain fundamental aspects concerning their structural and geometric properties remain to be fully explored.

Given an idempotent $Q\in\mathbb{B}(H)$, one may consider the set of operator distances from all projections to $Q$. A natural problem is to determine the minimum value, the maximum value, and the intermediate values of these distances. To this end,
a projection $m(Q)$, referred to as the matched projection of $Q$, is introduced for every idempotent $Q\in\mathbb{B}(H)$  in \cite{TXF02} as follows:
\begin{equation}\label{equ:m(q)}
m(Q)=\frac12\big(|Q^*|+Q^*\big)|Q^*|^\dag\big(|Q^*|+I\big)^{-1}\big(|Q^*|+Q\big),
\end{equation}
where $|Q^*|^\dag$ denotes the Moore-Penrose inverse of $|Q^*|$ \cite{XU,XS}. This can be explicitly written as
\begin{equation*}\label{formula for the MP inverse-01}|Q^*|^\dag=\left(P_{\mathcal{R}(Q)}P_{\mathcal{R}(Q^*)}P_{\mathcal{R}(Q)}\right)^\frac12.
\end{equation*}
Some investigations and applications of the matched projections can be found in \cite{TXF03,TXF02,TXF04,ZTX,ZFL}.

In this paper, we  investigate several applications of the matched projections. The primary application pertains to the operator distances from projections to a given idempotent.
Motivated by results obtained in \cite{TXF03} for matrices endowed with the Frobenius norm, the following problems are formulated:

\begin{problem}\label{problem-01}{\rm [Minimum value]} Let $Q\in\mathbb{B}(H)$ be an arbitrary non-projection idempotent. Is it true that
\begin{equation}\label{minimum distance-00}\|m(Q)-Q\|\le \|P-Q\|\end{equation}
for any projection $P\in\mathbb{B}(H)$?
\end{problem}

\begin{problem}\label{problem-02}{\rm [Maximum value]} Let $Q\in\mathbb{B}(H)$ be an arbitrary non-projection idempotent. Is it true that
\begin{equation*}\label{maximum distance-00}\|P-Q\|\le \|I-m(Q)-Q\|\end{equation*}
for any projection $P\in\mathbb{B}(H)$?
\end{problem}

\begin{problem}\label{prob:intermediate value}{\rm [Intermediate values]}   Let $Q\in\mathbb{B}(H)$ be an arbitrary non-projection idempotent, and let $\alpha\in
\big[\|m(Q)-Q\|,\|I-m(Q)-Q\|\big]$ be arbitrary. Is it true that there exists a projection $P\in \mathbb{B}(H)$ such that $\|P-Q\|=\alpha$?
\end{problem}

The positive answer to Problem~\ref{problem-02} is provided in \cite[Theorem~2.3]{ZTX}.
Partial results for Problem~\ref{problem-01} exist in the context of adjointable operators on Hilbert $C^*$-modules, where inequality \eqref{minimum distance-00}  holds when $(P,Q)$ forms a quasi-projection pair; see  \cite[Definition~2.2 and Theorem~4.15]{TXF02}.

Section~\ref{sec:distance} of this paper aims to resolve these problems comprehensively. Specifically, Theorem~\ref{thm:minimum distance} affirmatively answers  Problem~\ref{problem-01}, demonstrating that the distance from  $m(Q)$ to $Q$ achieves the minimum value among all operator distances from projections to $Q$.
 Furthermore, Theorem~\ref{thm:min to max} shows that the maximum distance  can be derived from the minimum distance. offering a novel proof for existing results; see Corollary~\ref{cor:norm of I-mq-q} for details.

Leveraging properties of the matched projection established in \cite{TXF02}, Theorem~\ref{thm:simplify R and N} positively addresses Problem~\ref{prob:intermediate value} by adapting strategies from finite-dimensional matrix analyses employed in
\cite[Theorem~4.5]{TXF03}.

Let $Q\in\mathbb{B}(H)$ be a non-projection idempotent.   There has been extensive research on the relationships between $Q$ and the pair $P_{\mathcal{R}(Q)}$ and $P_{\mathcal{N}(Q)}$, where $P_{\mathcal{R}(Q)}$ and $P_{\mathcal{N}(Q)}$ denote the projections from $H$ to the range and null space of $Q$, respectively. In particular, it is known that
the $C^*$-algebra generated by $Q$ and $I$  coincides with the $C^*$-algebra generated by $P_{\mathcal{R}(Q)}$ and $P_{\mathcal{N}(Q)}$; see Remark~\ref{rem:CstarQ&I} and \cite[Theorem~5.1]{BS}. Although the titles of references \cite{KRS} and \cite{RS} suggest the study of $C^*$-algebras generated by multiple idempotents and by two projections, respectively, they in fact refer to $C^*$-algebras generated by idempotents together with the identity operator, and by two projections together with the identity operator.

 To the best of our knowledge, research devoted to  the structure of the $C^*$-algebra generated by a single idempotent $Q$ on a Hilbert space, denoted by $C^*\{Q\}$, remains limited in the existing literature.  A  primary objective  of Section~\ref{subsec:structure} is to establish  a new representation formula for any idempotent $Q$, and subsequently provide a novel characterization of the  generated $C^*$-algebra $C^*\{Q\}$. This is achieved in Theorem~\ref{thm:new express of Q} by retaining the range projection $P_{\mathcal{R}(Q)}$ and replacing the null-space projection $P_{\mathcal{N}(Q)}$ with the matched projection $m(Q)$.

Recent work in \cite{TXF04} has established new block matrix representations for idempotents through the use of the matched projections.
Building on these representations, we provide a complete characterization of  all elements in $C^*\{Q\}$ by  using $4 \times 4$ block matrix forms in the case where $Q$ is a non-projection idempotent. This characterization is fully developed in Theorems~\ref{the:1 not in sigmaD} and \ref{th:1 in sigma}.

In Section~\ref{sec:universal-r-idempotents}, we  introduce  the universal $r$-idempotents for each parameter $r>1$. These idempotents exhibit a universal property, as formally stated in Definition~\ref{defn:universalr} and  alternatively characterized in Proposition~\ref{prop:alternativeup}. Notably, this universal property is distinct from the well-known property for pairs of projections described in \cite[Proposition~1.1]{RS} (see Remark~\ref{rem:differup}). For every $r>1$, an explicit construction of a universal $r$-idempotent  $Q_r$ is provided in Theorem~\ref{th:C*-homorphism2}.

Building upon the universal $r$-idempotent $Q_r$ and employing the $4\times 4$ block  matrix representations for all elements in $C^*\{Q\}$ where $Q$ is a non-projection idempotent $Q$,
Theorems~\ref{thm:nscondition4ui} and \ref{thm:4 r-idem nonblock} establish  necessary and sufficient conditions for $Q$ to be a universal $r$-idempotent.  Corollary~\ref{cor:Q&I-Quri} shows that for any non-projection $Q$, it is a universal $r$-idempotent if and only if $I-Q$ is also a universal $r$-idempotent. Moreover, it is proved in Theorem~\ref{non-uniqueness universal} that up to unitary equivalence, the universal $r$-idempotents are not unique for every $r>1$.

The elliptical range theorem is originally obtained in the matrix case, which is a powerful tool in dealing with the numerical ranges of matrices. The main purpose of Section~\ref{sub:operator version} is to set up an operator version of the elliptical range theorem. By considering various generalized numerical ranges in stead of the standard numerical range, several generalizations of the elliptical range theorem for quadratic operators have been established in  \cite{LPS,LPT}. Recall that an operator $A\in\mathbb{B}(H)$ is said to be quadratic if it satisfies
$$(A-a I_H)(A-b I_H)=0$$ for some complex numbers $a$ and $b$ \cite[Section 2.1]{WG}. Notable examples of quadratic operators include idempotents and square-zero operators. According to \cite[Theorem~1.1]{TW}, up to unitary equivalence, every quadratic operator can be represented in the block matrix form
$$aI\oplus b I\oplus \left(
                       \begin{array}{cc}
                         aI & D \\
                         0& bI\\
                       \end{array}
                     \right),$$
where $D$ is a positive operator with $\mathcal{N}(D)=\{0\}$. Note that all four entries of the above $2\times 2$ block matrix belong to the commutative $C^*$-algebra $C^*\{D,I\}\cong C\big(\sigma(D)\big)$. This motivates a new generalization of the elliptical range theorem by studying operators of the form
\begin{equation*}\label{equ:simple4fD}\begin{pmatrix}
        f_{11}(D) & f_{12}(D) \\
        f_{21}(D) & f_{22}(D)
       \end{pmatrix},\end{equation*}
where $D\in \mathbb{B}(H)_{\text{sa}}$ and $f_{ij}\in C\big(\sigma(D)\big)$ for $1\le i,j\le 2$. Earlier contributions to this aspect of the problem are scattered throughout the literature, including \cite{Klaja,Lenard}. A major obstacle in these studies is the diagonalizing of operator matrices. This process necessitates the construction of rather intricate two-parameter unitary matrices whose entries are Borel measurable functions (including continuous ones) of the first parameter (see \cite[Lemma~3.1]{Klaja} and \cite[Section~2]{Lenard}). In Section~\ref{sub:operator version}, we proceed by applying Lemma~\ref{lem:the norm of T+}  to overcome this technical difficulty, which subsequently leads to a general elliptical range theorem as presented in Theorem~\ref{thm:the numerical of matrix02}.

As an application, in Section~\ref{subsec:non-quadratic} we examine the numerical ranges of certain non-quadratic operators constructed from non-projection idempotents and their matched projections. The main goal is to  characterize their numerical ranges, thereby extending the results from the previous literatures. This is achieved by constructing a non-quadratic operator $T_Q$ in \eqref{equ:nonquadopt} and describing its numerical range in detail; see Theorems~\ref{thm:the num of Q2},  \ref{thm:biggerintT}, \ref{thm:1stcloseT} and Corollary~\ref{cor:num of Qr}. We also aim to identify meaningful examples and reveal new interesting phenomena. To this end, an operator $S_r$ is introduced as in \eqref{equ:Sr} for each parameter $r>1$. Theorems~\ref{thm:Srnoted} and \ref{thm:ncno} demonstrate that for this operator, the closure of its numerical range is not an elliptical disk, and the numerical range itself is neither closed nor open.

Throughout the remainder of this paper, $m(Q)$ denotes the matched projection of an idempotent $Q$. For any $C^*$-algebra $\mathfrak{B}$, $M_n(\mathfrak{B})$ denotes the $n\times n$ matrix algebra over $\mathfrak{B}$. For every subset $E$ of the complex field $\mathbb{C}$, let $\overline{E}$, $\text{int}(E)$, $\partial E$ and $\mbox{conv}(E)$ be the closure, interior,  boundary, and convex hull of $E$, respectively.

The rest of this paper is organized as follows. Section~\ref{sec:distance} investigates the operator distances from projections to a given idempotent. Section~\ref{sec: the min-maximum distances} provides positive answers to Problems~\ref{problem-01} and \ref{problem-02}. Section~\ref{sec:intermediate values} focuses on the intermediate values problem, offering a positive answer to Problem~\ref{prob:intermediate value}. Section~\ref{sec:C-algebras generated} is devoted to the study of the $C^*$-algebra generated by a single non-projection idempotent. A complete characterization of the elements of this $C^*$-algebra is presented in Section~\ref{subsec:structure}. Section~\ref{sec:universal-r-idempotents} introduces the notion of universal $r$-idempotents for each parameter $r>1$, providing a necessary and sufficient condition for a non-projection idempotent to be a universal $r$-idempotent. Section~\ref{sec:numerical range} develops a new characterization of the numerical range. This work leads to an operator version  of the elliptical range theorem in Section~\ref{sub:operator version}, which is then applied in Section~\ref{subsec:non-quadratic} to study the numerical ranges of certain non-quadratic operators constructed from non-projection idempotents and their matched projections.

\section{The operator distances from projections to an idempotent}\label{sec:distance}

\subsection{The minimum and the maximum operator distances}\label{sec: the min-maximum distances}

Suppose $Q\in \mathbb{B}(H)$ is a non-projection idempotent. Then $K_1 \neq \{0\}$ and $K_2 \neq \{0\}$, where
\begin{equation}\label{equ:K_1 and K_2}
K_1 = \mathcal{R}(Q), \quad K_2 = \mathcal{N}(Q^*).
\end{equation}
Let $U_{P_{K_1}}$ denote the unitary operator induced by the projection $P_{K_1}$ (see \eqref{equ:unitary operator induced by P}).
From \cite[Theorem~3.1]{TXF02}, we obtain
\begin{align}\label{blocked form of Q}&P_{K_1}= U_{P_{K_1}}^*\left(
\begin{array}{cc}
I_{K_1} & 0 \\
0 & 0 \\
\end{array}
\right)U_{P_{K_1}},\quad Q =U_{P_{K_1}}^*\left(
\begin{array}{cc}
I_{K_1} & A \\
0 & 0 \\
\end{array}
\right)U_{P_{K_1}},\\
\label{equ:blocked fromula for m q}&m(Q)=\frac12 U_{P_{K_1}}^*\left(\begin{array}{cc}
(B+I_{K_1})B^{-1} & B^{-1}A\\
A^*B^{-1} & A^*\big[B(B+I_{K_1})\big]^{-1}A \\
\end{array}\right) U_{P_{K_1}},
\end{align}
where $A \in \mathbb{B}(K_2, K_1)$ and $B \in \mathbb{B}(K_1)$ are defined by
\begin{equation}\label{equ:1st defn of AB}A=Q|_{K_2},\quad B=(I_{K_1}+AA^*)^\frac12.\end{equation}

To determine the minimum distance from projections to an idempotent $Q$, we require the norm formula for $m(Q)-Q$.
\begin{lemma}\label{lem:norm of mq-q}{\rm \cite[Theorem~4.17]{TXF02}}
 For any idempotent $Q\in \mathbb{B}(H)$,
\begin{equation}\label{norm of q minus m q}
\|m(Q) - Q\| = \frac{1}{2} \left[ \|Q\| - 1 + \sqrt{\|Q\|^2 - 1} \right].
\end{equation}
\end{lemma}

%%%%

Building upon the proof technique in \cite[Lemma~2.2]{ZTX}, we establish the key technical result of this section.
\begin{lemma}\label{lem:technical lem}Given non-zero Hilbert spaces $K_i$ $(i=1,2)$ and $A\in \mathbb{B}(K_2,K_1)$, let $T\in \mathbb{B}(K_1\oplus K_2)$ be defined by
\begin{equation}\label{this form of T-01}T=\left(
      \begin{array}{cc}
        I_{K_1} & A \\
        A^* & I_{K_2}+2A^*A \\
      \end{array}
    \right).\end{equation}
Then
\begin{equation}\label{norm of this T}\|T\|=1+\|A\|^2+\|A\|\sqrt{1+\|A\|^2}.
\end{equation}
\end{lemma}
\begin{proof}If $A=0$, the conclusion holds trivially. So, we may assume that  $A\ne 0$. Let $A=V_A|A|$ be the polar decomposition, where $V_A\in \mathbb{B}(K_2,K_1)$ is a partial isometry satisfying $\mathcal{R}(V_{A})=\overline{\mathcal{R}(A)}$ and $\mathcal{R}(V^*_{A})=\overline{\mathcal{R}(A^*)}$. Define
$$X=\left(
      \begin{array}{cc}
        I_{K_1} &  \\
         & V_A \\
      \end{array}
    \right)T\left(
      \begin{array}{cc}
        I_{K_1} &  \\
         & V_A^* \\
      \end{array}
    \right),\ Y=\left(
      \begin{array}{cc}
        I_{K_1} &  \\
         & V_A^* \\
      \end{array}
    \right)X\left(
      \begin{array}{cc}
        I_{K_1} &  \\
         & V_A \\
      \end{array}
    \right).$$
Direct computation yields
\begin{equation}\label{expressions of X and Y 2times2}X=\left(
      \begin{array}{cc}
        I_{K_1} & |A^*| \\
        |A^*| & V_AV_A^*+2 |A^*|^2 \\
      \end{array}
    \right),\quad Y=\left(
      \begin{array}{cc}
        I_{K_1} & A \\
        A^* & V_A^*V_A+2 A^*A \\
      \end{array}
    \right).\end{equation}

Let $\mathfrak{B}$ denote the $C^*$-subalgebra of $\mathbb{B}(K_1)$
generated by $I_{K_1}, |A^*|$ and $V_AV_A^*$.
Since
\begin{equation}\label{equ:prepared commutativity}V_AV_A^*|A^*|=|A^*|=|A^*|V_AV_A^*,\end{equation}
$\mathfrak{B}$ forms a unital commutative $C^*$-algebra. Consequently, there exist a compact Hausdorff space $\Omega$ and a unital $C^*$-algebra isomorphism $\rho:\mathfrak{B}\to
C(\Omega)$ \cite[Section~1.1]{Pedersen}. Let
\begin{equation}\label{prepare for a t b t}a=\rho(|A^*|),\quad  f=\rho(V_AV_A^*)=\chi_{\Omega_1},\end{equation}
where $\chi_{\Omega_1}$ denotes the characteristic function with $\Omega_1=\{t\in\Omega: f(t)=1\}$, which is both open and closed in $\Omega$.
Since every commutative $C^*$-algebra (as well as $M_2(\mathbb{C})$) is nuclear \cite[Theorem~6.3.9 and 6.4.15]{Murphy}, and
$$M_2(\mathfrak{B})\cong M_2(\mathbb{C})\otimes \mathfrak{B}\cong M_2(\mathbb{C})\otimes C(\Omega)\cong M_2\big(C(\Omega)\big),$$
the norm of the operator $X$ can be characterized as
\begin{equation*}\|X\|=\max_{t\in \Omega}\|X(t)\|=\max_{t\in \Omega}\left\|\left(
                                            \begin{array}{cc}
                                              X_{11}(t) & X_{12}(t) \\
                                              X_{12}(t) & X_{22}(t) \\
                                            \end{array}
                                          \right)\right\|,\end{equation*}
where $a\in C(\Omega)$ is defined by \eqref{prepare for a t b t}, and for each $t\in\Omega$ we have
\begin{align*}X_{11}(t)\equiv 1,\quad
X_{12}(t)=a(t),\quad
X_{22}(t)=\chi_{\Omega_1}(t)+2a^2(t).
\end{align*}

\textbf{Case 1:}\quad $V_AV_A^*=I_{K_1}$. Here $\Omega_1=\Omega$, so $X_{22}(t)=1+2a^2(t)$ for $t\in\Omega$. Let $\lambda_1(t)$ and $\lambda_2(t)$ denote the eigenvalues of $X(t)$. Direct computation yields
\begin{align*}\lambda_i(t)=1+a^2(t)\pm a(t)\sqrt{a^2(t)+1},\quad i=1,2.
\end{align*}
Since $X(t)$ is positive semi-definite and $a(t)\geq 0$, we obtain
\begin{align}\label{norm of X t}\|X(t)\|=&\max\{\lambda_1(t),\lambda_2(t)\}=1+a^2(t)+a(t)\sqrt{a^2(t)+1}.
\end{align}
Therefore,
\begin{align*}\|X\|=&\max_{t\in \Omega}\left[1+a^2(t)+a(t)\sqrt{a^2(t)+1}\right]\\
=&1+\|\rho(|A^*|)\|^2+\|\rho(|A^*|)\|\sqrt{\|\rho(|A^*|)\|^2+1}\\
=&1+\|A\|^2+\|A\|\sqrt{1+\|A\|^2}.
\end{align*}

\textbf{Case 2:}\quad $V_AV_A^*\ne I_{K_1}$. Here $\Omega_1\ne \Omega$. Set $\Omega_2=\Omega\setminus \Omega_1$. From \eqref{equ:prepared commutativity}, we have $\chi_{\Omega_1} (t)a(t)=a(t)$ for all $t\in\Omega$, where $a\in C(\Omega)$ is defined in \eqref{prepare for a t b t}. Consequently, $a(t)=0$ for all $t\in \Omega_2$. Thus for each $t\in \Omega_2$,
$$X(t)=\left(
         \begin{array}{cc}
           1 & 0 \\
           0 & 0 \\
         \end{array}
       \right),$$
yielding $\|X(t)\|=1$ for such $t$. Since $\chi_{\Omega_1}(t)=1$ for each $t\in\Omega_1$, as established in \textbf{Case 1}, $\|X(t)\|$ is given by \eqref{norm of X t} for $t\in\Omega_1$. Hence,
\begin{align*}
  \|X\|&=\max\big\{\max_{t\in\Omega_1}\|X(t)\|, \max_{t\in\Omega_2}\|X(t)\|\big\}=\max_{t\in\Omega_1}\|X(t)\| \\
   &=1+\|A\|^2+\|A\|\sqrt{1+\|A\|^2}.
\end{align*}
In both cases, we conclude
\begin{equation*}\|X\|=1+\|A\|^2+\|A\|\sqrt{1+\|A\|^2}.
\end{equation*}

It remains to demonstrate $\|T\|=\|X\|$. From the second equation in \eqref{expressions of X and Y 2times2} and the definitions of $X$ and $Y$, we have
$$1=\|I_{K_1}\|\le \|Y\|\le \|X\|\le \|T\|.$$ Moreover,  the expressions of $T$ and $Y$ given in \eqref{this form of T-01} and \eqref{expressions of X and Y 2times2}, respectively, yield
$$T=Y+P,\quad Y\ge 0,\quad P\ge 0, \quad YP=0,$$
where $P=\mbox{diag}\big(0, I_{K_2}-V_A^*V_A\big)$. This implies
$$\|T\|=\max\{\|Y\|,\|P\|\}=\|Y\|.$$
Thus $\|T\|=\|X\|$ as required.
\end{proof}

As a consequence of Lemma~\ref{lem:technical lem}, we obtain the following characterization of the minimum operator distance.
\begin{theorem}\label{thm:minimum distance} For any idempotent $Q\in\mathbb{B}(H)$,
$\|m(Q)-Q\|$ represents the minimum operator distance from $Q$ to all projections on $H$.
\end{theorem}
\begin{proof}For any projection $P\in\mathbb{B}(H)$, we have $m(P)=P$. Thus, the conclusion
holds trivially when $Q$ is a projection. We now consider the case where $Q\in\mathbb{B}(H)$ is a non-projection idempotent, decomposed as in \eqref{blocked form of Q}, with $K_1$ and $K_2$ defined by \eqref{equ:K_1 and K_2}.

For any projection $P\in\mathbb{B}(H)$, define
\begin{equation}\label{equ:defn of S P Q}S_{Q,P}=(P-Q)^*(P-Q)+(I-P-Q)^*(I-P-Q).\end{equation}
Direct computation yields
\begin{equation}\label{equ: S Q P}S_{Q,P}=I-Q-Q^*+2Q^*Q,\end{equation}
which shows that $S_{Q,P}$ is independent of the choice of $P$. From \eqref{equ: S Q P} and \eqref{blocked form of Q}, we obtain
$$S_{Q,P}=U_{K_1}^*TU_{K_1},$$
where $T$ is given by \eqref{this form of T-01}. Applying \eqref{norm of this T} gives
\begin{equation}\label{norm of S Q P}\|S_{Q,P}\|=1+a^2+a\sqrt{1+a^2},
\end{equation}
where $a=\|A\|$. Meanwhile, observe that
\begin{equation*}\|I-P-Q\|\le \|I-2P\|+\|P-Q\|= 1+\|P-Q\|.\end{equation*}
From \eqref{equ:defn of S P Q}, we deduce
\begin{align}\label{equ:small technique-1}\|S_{Q,P}\|\le& \|P-Q\|^2+\|I-P-Q\|^2\\
\label{equ:small technique-2}\le& \|P-Q\|^2+\big(1+\|P-Q\|\big)^2.
\end{align}
Combining \eqref{norm of S Q P}  with \eqref{equ:small technique-2} yields
$$a^2+a\sqrt{1+a^2}\le 2\|P-Q\|^2+2\|P-Q\|.$$
Consequently,
$$\big(\sqrt{1+a^2}+a\big)^2=1+2\left[a^2+a\sqrt{1+a^2}\right]\le \big(2\|P-Q\|+1\big)^2,$$
and therefore
$$\sqrt{1+a^2}+a\le 2\|P-Q\|+1.$$
Note that $\|Q\|=\sqrt{1+\|A\|^2}=\sqrt{1+a^2}$. Using \eqref{norm of q minus m q}, we obtain
\begin{equation*}\|m(Q)-Q\|=\frac{\sqrt{1+a^2}+a-1}{2}\le \|P-Q\|.
\end{equation*}
This establishes that $\|m(Q)-Q\|$ indeed gives the minimum operator distance.
\end{proof}

Our next theorem demonstrates that the maximum operator distance can be derived from the minimum operator distance.
\begin{theorem}\label{thm:min to max} Let $Q\in \mathbb{B}(H)$ be an idempotent and $P_0\in\mathbb{B}(H)$ be a projection such that
$\|P_0-Q\|$ achieves the minimum operator distance from $Q$ to all projections on $H$. Then
\begin{equation}\label{equ:min to max}\|I-P_0-Q\|=1+\|P_0-Q\|,
\end{equation}
which is the maximum operator distance from $Q$ to all projections on $H$.
\end{theorem}
\begin{proof}From the Krein-Krasnoselskii-Milman equality (see \cite[Lemma~4.1]{XY2}), we have
$\|P_1-P_2\|\le 1$
for any projections $P_1$ and $P_2$ on $H$. Thus, we may assume $Q$ is a non-projection idempotent.

Define $S_{Q,P}$ as in \eqref{equ:defn of S P Q} for any projection $P$ on $H$.
By Theorem~\ref{thm:minimum distance}, \eqref{norm of q minus m q}, and letting $a = \|A\|$, we obtain \begin{align*}\|P_0-Q\|=\frac{\sqrt{1+a^2}-1+a}{2},
\end{align*}
which implies
\begin{equation}\label{norm of 001}\|P_0-Q\|^2+\left(1+\|P_0-Q\|\right)^2=1+a^2+a\sqrt{1+a^2}\end{equation}
through direct computation. Combining \eqref{equ:small technique-1}, \eqref{equ:small technique-2}, \eqref{norm of 001}, and \eqref{norm of S Q P} yields
\begin{align*}\|S_{Q,P_0}\|\le& \|P_0-Q\|^2+\|I-P_0-Q\|^2\\
\le& \|P_0-Q\|^2+\left(1+\|P_0-Q\|\right)^2=\|S_{Q,P_0}\|,
\end{align*}
establishing the validity of \eqref{equ:min to max}.

For any projection $P\in\mathbb{B}(H)$, we have $\|P-P_0\|\le 1$, so
$$\|P-Q\|\le \|P-P_0\|+\|P_0-Q\|\le 1+\|P_0-Q\|=\|I-P_0-Q\|.$$
This confirms that $\|I-P_0-Q\|$ indeed represents the maximum distance.
\end{proof}

Combining Theorems~\ref{thm:minimum distance} and \ref{thm:min to max}, we immediately recover the main results of \cite{ZTX}:
\begin{corollary}\label{cor:norm of I-mq-q}{\rm \cite[Lemma~2.2 and Theorem~2.3]{ZTX}} Let $Q\in \mathbb{B}(H)$ be an idempotent. Then
\begin{equation*}\label{norm of max wrt min}\|I-m(Q)-Q\|=1+\|m(Q)-Q\|,
\end{equation*}
which is the maximum operator distance from $Q$ to all projections on $H$.
\end{corollary}

\begin{remark} It should be noted that the converse of Theorem~\ref{thm:min to max}  does not hold. Consider the following example: there exist an idempotent $Q$ and a projection $P$ such that $\|I-P-Q\|$ achieves the maximum operator distance from $Q$ to all projections, but $\|P-Q\|$ is not the minimum operator distance from $Q$ to all projections.
\end{remark}

\begin{example} For any $a\in (0, \frac43)$, let $A\in M_2(\mathbb{C})$ and $Q\in M_4(\mathbb{C})$ be defined by
$$A=\left(
    \begin{array}{cc}
      a & 0 \\
      0 & 0 \\
    \end{array}
  \right),\quad Q=\left(
                    \begin{array}{cc}
                      I_2 & A \\
                      0 & 0 \\
                    \end{array}
                  \right).$$
Since $0<\|A\|=a<\frac43$, by \eqref{norm of q minus m q} and Corollary~\ref{cor:norm of I-mq-q}, we have
\begin{equation*}\|m(Q)-Q\|=\frac{\sqrt{1+a^2}+a-1}{2}<1,\quad \|I_4-m(Q)-Q\|=1+\|m(Q)-Q\|>1.\end{equation*}
Define $P_0,P,T\in M_4(\mathbb{C})$ by
$$P_0=\mbox{diag}(0,0,0,1),\quad  P=m(Q)+P_0,\quad T=I_4-P-Q.$$
Observe that $Q$ can also be partitioned  as
$$Q=\left(
      \begin{array}{ccc|c}
        1 & 0 & a & 0 \\
        0 & 1 & 0 & 0 \\
        0 & 0 & 0 & 0 \\\hline
        0 & 0 & 0 & 0 \\
      \end{array}
    \right)=\left(
              \begin{array}{cc}
                Q_1 & 0 \\
                0 & 0 \\
              \end{array}
            \right),$$
which gives $m(Q)=m(Q_1)\oplus 0$. Hence,
$P_0Q=QP_0=0_4$, and $P_0m(Q)=m(Q)P_0=0_4$. It follows that $P$ is a projection, and
 \begin{align*}\|P-Q\|^2=&\|(P-Q)^*(P-Q)\|=\big\|\big(m(Q)-Q\big)^*\big(m(Q)-Q\big)+P_0\big\|\\
 =&\max\big\{ \|m(Q)-Q\|^2, \|P_0\|\big\}=\max\big\{ \|m(Q)-Q\|^2, 1\big\}=1.
 \end{align*}
 Thus,
 $$\|P-Q\|=1>\|m(Q)-Q\|,$$
 proving that $\|P-Q\|$ is not the minimum operator distance.

Next, note that $T=(I_4-P_0)-m(Q)-Q$, so $P_0T=TP_0=0_4$. Therefore, similar reasoning gives
$$\|I_4-m(Q)-Q\|=\|T+P_0\|=\max\{\|T\|,\|P_0\|\}=\max\{\|T\|,1\}.$$
Since $\|I_4-m(Q)-Q\|>1$, it follows that $\|T\|=\|I_4-m(Q)-Q\|$. Hence, $\|I_4-P-m(Q)\|$ achieves the maximum operator distance.
\end{example}

\begin{remark}
By the same reason as in \cite[Remark~2.1]{ZTX}, Theorems~\ref{thm:minimum distance} and \ref{thm:min to max} hold for an adjointable idempotent and projections on a Hilbert $C^*$-module.
\end{remark}

\subsection{The intermediate values of the operator distances}\label{sec:intermediate values}

For an idempotent $Q\in\mathbb{B}(H)$, define $\mbox{ran}_Q$  by
\begin{equation}\label{defn of ran Q}
\mbox{ran}_Q=\big\{\|P-Q\|: \mbox{$P$ is a projection on $H$}\big\}.
\end{equation}
By Theorem~\ref{thm:minimum distance} and Corollary~\ref{cor:norm of I-mq-q}, we have $\mbox{ran}_Q\subseteq [\mbox{min}_Q, \mbox{max}_Q]$, where
\begin{equation*}\mbox{min}_Q=\|m(Q)-Q\|,\quad \mbox{max}_Q=\|I-m(Q)-Q\|.
\end{equation*}
A natural question is whether the reverse inclusion $[\mbox{min}_Q, \mbox{max}_Q]\subseteq \mbox{ran}_Q$ holds.
This inclusion can fail when $Q$ is a projection. For instance, if $Q$ is a trivial projection, then
 $[\mbox{min}_Q, \mbox{max}_Q]=[0,1]$, while the norm  $\|P-Q\|$ takes only the values 0 or 1 for any projection $P$. Hence, we restrict our attention to non-projection idempotents.

\begin{lemma}\label{lem:two relationships with m q}{\rm \cite[Theorems~3.7 and 3.14]{TXF02}} Let $Q\in \mathbb{B}(H)$ be an idempotent. Then \begin{equation*}\label{2 equs wrt m q}m(Q^*)=m(Q), \quad m(I-Q)=I-m(Q).\end{equation*}
\end{lemma}

%%%

\begin{lemma}\label{lem:partial intermediate} For any idempotent $Q\in\mathbb{B}(H)$, define
\begin{equation}\label{equ: defn of lambda q m q}\lambda_Q=\|P_{\mathcal{R}(Q)}-Q\|,\quad \mu_Q=\|I-P_{\mathcal{R}(Q)}-Q\|.
\end{equation}
Then $\lambda_Q<\mu_Q$ and
\begin{equation*}\left[\mbox{min}_Q,\lambda_Q\right]\cup \left[\mu_Q, \mbox{max}_Q\right]\subseteq \mbox{ran}_Q.
\end{equation*}
\end{lemma}
\begin{proof}We first establish $\lambda_Q<\mu_Q$. If $Q$ is a projection, then $\lambda_Q=0$ and $\mu_Q=1$, confirming the inequality.
Now suppose $Q\in\mathbb{B}(H)$ is a non-projection idempotent. Let $K_1$ and $K_2$ be defined by \eqref{equ:K_1 and K_2}, with $P_{K_1}$ and $Q$ partitioned as in \eqref{blocked form of Q}, where $A$ is given by \eqref{equ:1st defn of AB}. Then
\begin{align*}\lambda_Q=\|A\|<\sqrt{1+\|A\|^2}\le \left\|\left(
                                                           \begin{array}{cc}
                                                             -I_{K_1} & -A \\
                                                             0 & I_{K_2} \\
                                                           \end{array}
                                                         \right)\right\|=\mu_Q.
\end{align*}

Next, we demonstrate $\left[\mbox{min}_Q,\lambda_Q\right]\subseteq \mbox{ran}_Q$. This holds trivially when $Q$ is a projection. Suppose that $Q$ is a non-projection idempotent. Note that $P_{K_1}$ and $Q$ are homotopy equivalent as idempotents (see the proof of \cite[Lemma~4.4]{TXF03}). Similarly, $Q$ and $m(Q)$ are homotopy equivalent by \cite[Theorem~3.1]{TXF02}. Thus, there exists a norm-continuous path $Q_t$ ($t\in [0,1]$) of idempotents in $\mathbb{B}(H)$ connecting $m(Q)$ to $P_{K_1}$. Define
$$f(t)=\left\|m(Q_t)-Q\right\|,\quad  t\in [0,1].$$
By \cite[Remark~3.5]{TXF02}, the mapping from an idempotent to its matched projection is norm-continuous. Therefore, $f$ is continuous on $[0,1]$, and since $f(0)=\mbox{min}_Q$, $f(1)=\lambda_Q$, and $m(Q_t)$ is a projection for each $t\in [0,1]$, we conclude $\left[\mbox{min}_Q,\lambda_Q\right]\subseteq \mbox{ran}_Q$.

Finally, we verify $\left[\mu_Q, \mbox{max}_Q\right]\subseteq \mbox{ran}_Q$. From Lemma~\ref{lem:two relationships with m q}, we have
$$\|I-m(Q)-Q\| =\|m(I-Q^*)-Q\|.$$
Observe that
$$P_{\mathcal{R}(I-Q^*)}=P_{K_2}=I-P_{K_1},$$
so $I-P_{K_1}$ and $m(I-Q^*)$ are homotopy equivalent. Applying analogous reasoning shows that $\left[\mu_Q, \mbox{max}_Q\right]\subseteq \mbox{ran}_Q$.
\end{proof}

\begin{lemma}\label{lem:infinite-dimensional case}For any non-projection idempotent $Q\in \mathbb{B}(H)$, define $\mbox{ran}_Q$ by \eqref{defn of ran Q}, and $\lambda_Q$, $\mu_Q$  by \eqref{equ: defn of lambda q m q}. Then $\left[\lambda_Q,\mu_Q\right)\subseteq \mbox{ran}_Q$.
\end{lemma}
\begin{proof}As before, let $K_1$ and $K_2$ be defined by \eqref{equ:K_1 and K_2}. For any projection $P\in \mathbb{B}(H)$, we observe that
$$\|I-P-(I-Q^*)\|=\|P-Q^*\|=\|(P-Q^*)^*\|=\|P-Q\|,$$
which implies $\mbox{ran}_{{}_{I-Q^*}}=\mbox{ran}_Q$. Furthermore,
\begin{align*}&\lambda_{{}_{I-Q^*}}=\|P_{\mathcal{R}(I-Q^*)}-(I-Q^*)\|=\|P_{K_1}-Q^*\|=\lambda_Q,\\
&\mu_{{}_{I-Q^*}}= \|I-P_{\mathcal{R}(I-Q^*)}-(I-Q^*)\|=\|I-P_{K_1}-Q^*\|=\mu_Q.
\end{align*}
By replacing $Q$ with $I-Q^*$ if necessary, we may assume $|K_1|\le |K_2|$, where $|K_1|$ and $|K_2|$ denote the cardinalities of $K_1$ and $K_2$, respectively. Thus, up to unitary equivalence, $K_2=K_1\oplus K_3$ for some Hilbert space $K_3$. Consequently, $H=K_1\oplus K_1\oplus K_3$, with $P_{K_1}$ and $Q$ given by
$$P_{K_1}=\begin{pmatrix} I_{K_1} &  &  \\ & 0 & \\ & & 0 \end{pmatrix},\quad Q=\begin{pmatrix} I_{K_1} & A_1 & A_2 \\ 0 & 0 & 0 \\ 0 & 0 & 0 \end{pmatrix}=\left(\begin{array}{c|c} I_{K_1} & A \\ \hline 0 & 0 \end{array}\right),$$
where $A=(A_1,A_2)$ with $A_1\in\mathbb{B}(K_1)$ and $A_2\in \mathbb{B}(K_3, K_1)$.

For each $t\in [0,1]$, define the projection $P_t\in\mathbb{B}(K_1\oplus K_1)$ by
$$P_t=\begin{pmatrix} \cos^2\left(\frac{\pi t}{2}\right) I_{K_1} & \sin\left(\frac{\pi t}{2}\right) \cos\left(\frac{\pi t}{2}\right) I_{K_1} \\ \sin\left(\frac{\pi t}{2}\right) \cos\left(\frac{\pi t}{2}\right) I_{K_1} & \sin^2\left(\frac{\pi t}{2}\right) I_{K_1} \end{pmatrix}.$$

\textbf{Case 1:}\quad $K_3\ne 0$. Choose  an arbitrary orthonormal basis $\{e_j: j\in J\}$ for $H_3$, and let $\mathcal{D}$ consist of all finite subsets of $J$ (including the empty set $\emptyset$). Define $P^{(\emptyset)}=0$, and for each non-empty $S\in\mathcal{D}$, let $P^{(S)}$ be the projection onto $\mbox{span}\{e_j:j\in S\}$. For $t\in [0,1]$ and $S\in\mathcal{D}$, define
\begin{equation}\label{defn of f S t}P_{t,S}=\begin{pmatrix} P_t &0 \\0 & P^{(S)} \end{pmatrix},\quad f_S(t)=\|P_{t,S}-Q\|.\end{equation}
Notably, every $P_{t,S}$ is a projection in $\mathbb{B}(H)$, and when $S$ is fixed, the function $f_S$ is continuous on $[0,1]$. Thus,
\begin{equation}\label{for contained in the range}\mathcal{R}(f_S)=\big\{f_S(t):0\le t\le 1\big\}\subseteq \mbox{ran}_Q,\quad  S\in\mathcal{D}.\end{equation}

For each $S\in\mathcal{D}$, by the definition of $P_t$, we have
\begin{align}&\label{P 0  S}P_{0,S}-Q=\begin{pmatrix} 0 & -A_1 & -A_2\\ 0& 0 & 0 \\0 & 0 & P^{(S)} \end{pmatrix},\\
&\label{P 1 S}P_{1,S}-Q=\begin{pmatrix} -I_{K_1} & -A_1 & -A_2\\0 & I_{K_1} & 0 \\0 & 0 & P^{(S)} \end{pmatrix},
\end{align}
which yield $P_{0,\emptyset}-Q=P_{K_1}-Q$ and
$$P_{0,S}-Q=T_1^{(S)}(I-Q),\quad P_{1,S}-Q=T_2^{(S)}(I-P_{K_1}-Q),$$
where
$$T_i^{(S)}=\begin{pmatrix} I_{K_1} & & \\ & (i-1)I_{K_1} & \\ & & P^{(S)} \end{pmatrix},\quad i=1,2.$$
For all $S\in\mathcal{D}$ and $i=1,2$, it is clear that $\big\|T_i^{(S)}\big\|=1$. Hence, for any $S,S'\in\mathcal{D}$, we have
\begin{align*}f_S(0)&=\|P_{0,S}-Q\|=\big\|T_1^{(S)}(I-Q)\big\|\le \|I-Q\|=\|-Q\|\\
&=\|P_{K_1}(P_{1,S'}-Q)\|\le \|P_{1,S'}-Q\|=f_{S'}(1)\\
&=\big\|T_2^{(S')}(I-P_{K_1}-Q)\big\|\le \|I-P_{K_1}-Q\|=\mu_Q.
\end{align*}
In particular, $f_S(0)\le f_S(1)$ for all $S\in\mathcal{D}$.  Thus, by \eqref{for contained in the range} we conclude that  $$\big[f_S(0),f_S(1)\big]\subseteq \mathcal{R}(f_S).$$

For $S\subseteq S'$, we have $P^{(S)}P^{(S')}=P^{(S)}$. From \eqref{P 0  S} and \eqref{P 1 S},
\begin{equation*}P_{0,S}-Q=T_2^{(S)}\left(P_{0,S'}-Q\right), \quad P_{1,S}-Q=T_2^{(S)}\left(P_{1,S'}-Q\right).\end{equation*}
This means that $\{f_S(0)\}_{S\in\mathcal{D}}$ and $\{f_S(1)\}_{S\in\mathcal{D}}$ are increasing nets. Since $\lim\limits_{S\in\mathcal{D}}P^{(S)}=I_{H_3}$ in the strong operator topology, \eqref{P 1 S} implies
\begin{align*}\lim_{S\in\mathcal{D}} \big\|\big(P_{1,S}-Q\big)x-(I-P_{K_1}-Q)x\big\|=0
\end{align*}
for all $x\in H$, yielding
$$\mu_Q=\|I-P_{K_1}-Q\|\le \sup_{S\in\mathcal{D}}\big\|P_{1,S}-Q\big\|=\lim_{S\in \mathcal{D}}f_S(1)\le \mu_Q.$$
Therefore, $\mu_Q=\lim\limits_{S\in \mathcal{D}}f_S(1)$.

From \eqref{defn of f S t} and \eqref{P 0  S}, $f_{\emptyset}(0)=\|P_{0,\emptyset}-Q\|=\|P_{K_1}-Q\|=\lambda_Q$. Since $\{f_S(0)\}_{S\in\mathcal{D}}$ and $\{f_S(1)\}_{S\in\mathcal{D}}$ are increasing nets satisfying
$$f_{\emptyset}(0)=\lambda_Q,\quad \lim_{S\in \mathcal{D}}f_S(1)=\mu_Q,\quad f_S(0)\le f_{S'}(1), \quad S,S'\in\mathcal{D},$$
we have
\begin{equation*}[\lambda_Q,\mu_Q)\subseteq \bigcup_{S\in \mathcal{D}}\big[f_S(0),f_S(1)\big]\subseteq\bigcup_{S\in \mathcal{D}}\mathcal{R}(f_S)\subseteq \mbox{ran}_Q.\end{equation*}

\textbf{Case 2:}\quad $K_3=0$. In this case, up to unitary equivalence, we have
$$P_{K_1}=\left(
            \begin{array}{cc}
              I_{K_1} & 0 \\
              0 & 0 \\
            \end{array}
          \right),\quad Q=\left(
                            \begin{array}{cc}
                              I_{K_1} & A_1 \\
                              0 & 0 \\
                            \end{array}
                          \right).$$
For each $t\in [0,1]$, let  $P_t\in\mathbb{B}(K_1\oplus K_1)$ be defined as above. Define
$$ f(t)=\|P_{t}-Q\|,\quad t\in [0,1].$$
Then
$f(0)=\|P_{K_1}-Q\|=\lambda_Q$, and $f(1)=\|I-P_{K_1}-Q\|=\mu_Q$. Hence,
$$[\lambda_Q,\mu_Q]=[f(0),f(1)]\subseteq \mathcal{R}(f)\subseteq \mbox{ran}_Q.$$
This completes the proof.
\end{proof}

\begin{theorem}\label{thm:simplify R and N} The answer to Problem~\ref{prob:intermediate value} is affirmative.
\end{theorem}
\begin{proof}The conclusion follows directly from Lemmas~\ref{lem:partial intermediate} and \ref{lem:infinite-dimensional case}, which together establish that every value in the interval $\left[\mbox{min}_Q,\mbox{max}_Q\right]$ is attained as the operator distance between some projection $P$ and the non-projection idempotent $Q$.
\end{proof}

\section{$C^*$-algebras generated by non-projection idempotents}\label{sec:C-algebras generated}

Functional calculus for elements in a commutative $C^*$-algebra will be frequently used  in the remainder of this paper. To facilitate this, we introduce a necessary function. Define
\begin{equation*}\label{simplifited spectrum restriction}(0,1)^c=(-\infty,0]\cup [1,+\infty),
\end{equation*}
and let $\ell$ be the real-valued continuous  functions  defined by
\begin{equation}\label{defn of lt}
\ell(t)=\sqrt{t^2-t},\quad  t\in (0,1)^c.
\end{equation}
By definition, if $\mathfrak{A}$ is a  $C^*$-algebra and $D\in \mathfrak{A}_{\mbox{sa}}$ satisfying $\sigma(D)\subseteq (0,1)^c$,  then
\begin{align*}\ell(D)=\big(D(D-I)\big)^{\frac12}.
\end{align*}

\subsection{The structure of the  $C^*$-algebra generated by a single non-projection idempotent}\label{subsec:structure}

 We begin this subsection with the following definition.

\begin{definition}
For any $C^{*}$-algebra $\mathfrak{A}$ and its  subset $\mathfrak{F}$, let $C^{*}\mathfrak{F}$ denote the $C^{*}$-subalgebra of $\mathfrak{A}$ generated by $\mathfrak{F}$.
\end{definition}
For  an  idempotent  $Q\in \mathbb{B}(H)$,  the notation $C^*\{Q\}$ refers to the $C^*$-subalgebra of $\mathbb{B}(H)$ generated by $Q$ alone.  It should be noted that in \cite[Section~5]{BS}, the notation $C^*(Q)$ is used to denote the $C^*$-algebra generated by $Q$ together with the identity operator $I$. Thus,  $C^*(Q)=C^*\{Q,I\}$, which is distinct from $C^*\{Q\}$.

Given $T \in \mathbb{B}(H,K)$, recall that its Moore-Penrose inverse (abbreviated as M-P inverse) is  the unique element $X \in \mathbb{B}(K, H)$  satisfying
\begin{equation}\label{equ:4penrose}
  TXT=T,\quad XTX=X,\quad (TX)^*=TX, \quad (XT)^*=XT.
\end{equation}
It is known that $T$ has the M-P inverse if and only if $\mathcal{R}(T)$ is closed in $K$ \cite[Theorem~2.2]{XS}. In this case, $T$ is said to be M-P invertible and its M-P inverse is denoted by $T^\dag$.

\begin{lemma}\label{lem: A dag contained sub alg}\label{lem:mpincluded} Suppose $T\in\mathbb{B}(H)$ is M-P invertible. Then $T^\dag\in C^*\{T\}$.
\end{lemma}
\begin{proof} For a proof, see \cite[Corollary~4.22]{BS02} and \cite[Lemma~4.2]{XU}.
\end{proof}

Based on the above, we present the following useful lemma.

\begin{lemma}\label{lem:keepMP}Let $\mathfrak{A}$ and $\mathfrak{B}$ be $C^*$-subalgebras of $\mathbb{B}(H)$ and $\mathbb{B}(K)$, respectively, and let $\pi:\mathfrak{A}\to \mathfrak{B}$ be a $C^*$-algebra homomorphism. If $T\in \mathfrak{A}$ is  M-P invertible, then $\pi(T)$ is also M-P invertible in $\mathfrak{B}$ such that $\pi(T)^\dag=\pi(T^\dag)$.
\end{lemma}
\begin{proof}By Lemma~\ref{lem: A dag contained sub alg}, $T^\dag\in C^*\{T\}\subseteq \mathfrak{A}$, so $\pi(T^\dag)$ is well-defined. Applying the homomorphism $\pi$ to both sides of each of the four Penrose equations (see \eqref{equ:4penrose}) shows that $\pi(T^\dag)$ satisfies the defining conditions for the M-P inverse of $\pi(T)$. Hence, $\pi(T)^\dag=\pi(T^\dag)$.
\end{proof}

The following theorem establishes a new formula for any idempotent $Q$, and characterizes the structure of the generated $C^*$-algebra $C^*\{Q\}$.
\begin{theorem}\label{thm:new express of Q}
  Let $Q\in \mathbb{B}(H)$  be an idempotent. Then
    \begin{align}\label{equ:new express of Q}&Q=\big[P_{\mathcal{R}(Q)}\big(2m(Q)-I\big)P_{\mathcal{R}(Q)}\big]^\dag P_{\mathcal{R}(Q)}\big(2m(Q)-I\big),\\
\label{equ:CstarQ}&C^*\{Q\}=C^*\big\{m(Q),P_{\mathcal{R}(Q)}\big\}=C^*\big\{m(Q),I-P_{\mathcal{N}(Q)}\big\}.
  \end{align}
 \end{theorem}
  \begin{proof}
If $Q$ is a projection, then $P_{\mathcal{R}(Q)}=m(Q)=Q=Q^\dag$ and $P_{\mathcal{N}(Q)}=I-Q$. Thus,  both \eqref{equ:new express of Q} and \eqref{equ:CstarQ}
are trivially satisfied.

Now suppose  $Q$ is a non-projection idempotent.
Let $K_1$, $K_2$ and $U_{P_{K_1}}$ be defined as in \eqref{equ:K_1 and K_2} and \eqref{equ:unitary operator induced by P}, respectively. Define
\begin{align*}&T=2m(Q)-I,\quad S_1=P_{K_1} T P_{K_1},\quad S_2=P_{K_1} T.
\end{align*}
By \eqref{equ:blocked fromula for m q},  we have
\begin{align*}
  &U_{P_{K_1}}T U_{P_{K_1}}^* = \left(
                                                  \begin{array}{cc}
                                                    B^{-1} & B^{-1}A \\
                                                    A^*B^{-1} & A^*\big[B(B+I_{K_1})\big]^{-1}A-I_{K_2}\\
                                                  \end{array}
                                                \right),
\end{align*}
where  $A$ and $B$ are given by  \eqref{equ:1st defn of AB}. It follows that
\begin{align*}U_{P_{K_1}}S_2U_{P_{K_1}}^*=&U_{P_{K_1}}P_{K_1}U_{P_{K_1}}^*\cdot U_{P_{K_1}}TU_{P_{K_1}}^*
=\left(
                                                  \begin{array}{cc}
                                                    B^{-1} & B^{-1}A \\
                                                    0 & 0\\
                                                  \end{array}
                                                \right),\\
U_{P_{K_1}}S_1^\dag U_{P_{K_1}}^*=&\big(U_{P_{K_1}}S_1U_{P_{K_1}}^*\big)^\dag=\left(
                                                  \begin{array}{cc}
                                                    B^{-1} & 0 \\
                                                    0 & 0 \\
                                                  \end{array}
                                                \right)^\dag=\left(
                                                  \begin{array}{cc}
                                                    B & 0 \\
                                                    0 & 0 \\
                                                  \end{array}
                                                \right).
\end{align*}
Therefore, by
 \eqref{blocked form of Q},
\begin{equation*}
   U_{P_{K_1}}S_1^\dag S_2U_{P_{K_1}}^*=\left(
                                    \begin{array}{cc}
                                      I_{K_1} & A \\
                                      0 & 0 \\
                                    \end{array}
                                  \right)=
  U_{P_{K_1}}QU_{P_{K_1}}^*,
\end{equation*}
which implies $Q=S_1^\dag S_2$. Hence, \eqref{equ:new express of Q} holds.

Since
$S_1=2P_{K_1} m(Q) P_{K_1}-P_{K_1}$, and $S_2=2P_{K_1} m(Q)-P_{K_1}$,
we have
$$C^*\{S_1\}\subseteq  C^*\big\{m(Q),P_{K_1}\},\quad  S_2\in C^*\big\{m(Q),P_{K_1}\}.$$
From $Q=S_1^\dag S_2$ and  Lemma~\ref{lem:mpincluded}, it follows that
$Q\in  C^*\big\{m(Q),P_{K_1}\}$. Thus,
$$C^*\{Q\}\subseteq  C^*\big\{m(Q),P_{K_1}\}.$$

On the other hand, by \eqref{equ:m(q)} and  Lemma~\ref{lem:mpincluded},
\begin{align*}&m(Q) \in C^{*}\{Q\} \cdot C^{*}\{Q, I\} \cdot C^{*}\{Q\} = C^{*}\{Q\}.
 \end{align*}
Moreover, as shown in \cite[Theorem~1.3]{Koliha}, the projection  $P_{\mathcal{R}(Q)}$ can be expressed as
  \begin{equation}\label{equ:projection Prq}
    P_{\mathcal{R}(Q)}=Q(Q+Q^*-I)^{-1},
  \end{equation}
which implies  $$ P_{\mathcal{R}(Q)}\in  C^*\{Q\}\cdot C^*\{Q,I\}=C^*\{Q\}. $$
Therefore,
$$C^*\big\{m(Q),P_{\mathcal{R}(Q)}\big\}\subseteq C^*\{Q\}.$$
This proves the first equality in \eqref{equ:CstarQ}.

Since $C^*\{Q\}=C^*\{Q^*\}$, $P_{\mathcal{R}(Q^*)}=I-P_{\mathcal{N}(Q)}$, and by Lemma~\ref{lem:two relationships with m q} we have $m(Q^*)=m(Q)$, the second equality in \eqref{equ:CstarQ} follows directly from the first.
\end{proof}

\begin{remark}Let $Q\in\mathbb{B}(H)$ be an idempotent. From the second equation in \eqref{equ:CstarQ}, it is immediate that
$I\in C^*\{Q\}$ if and only if $P_{\mathcal{N}(Q)}\in C^*\{Q\}$.
\end{remark}

As a direct consequence of Theorem~\ref{thm:new express of Q}, we obtain the following result from \cite{KRS}.
\begin{corollary}{\rm \cite[Theorem~6]{KRS}} For every  idempotent $Q\in\mathbb{B}(H)$, there exist two projections $P_1,P_2\in \mathbb{B}(H)$ such that
$C^*\{Q,I\}=C^*\{P_1,P_2,I\}$.
\end{corollary}
\begin{proof}The conclusion  follows by taking $P_1=m(Q)$ and $P_2=P_{\mathcal{R}(Q)}$.
\end{proof}

\begin{remark}\label{rem:CstarQ&I}
Suppose  that $Q\in \mathbb{B}(H)$ is an idempotent. Since $\mathcal{R}(I-Q)=\mathcal{N}(Q)$, it follows from \eqref{equ:projection Prq} that
\begin{equation*}P_{\mathcal{N}(Q)}=(Q-I)(Q+Q^*-I)^{-1}.\end{equation*}
Hence, $P_{\mathcal{N}(Q)}\in C^*\{Q,I\}$, and therefore
$C^*\{P_{\mathcal{R}(Q)},P_{\mathcal{N}(Q)}\}\subseteq C^*\{Q,I\}$.
On the other hand, by \cite[Theorem~3.8]{Ando02}, we have $$Q=P_{\mathcal{R}(Q)}\big(P_{\mathcal{R}(Q)}+ P_{\mathcal{N}(Q)}\big)^{-1},$$
which implies
$I,Q\in C^*\{P_{\mathcal{R}(Q)},P_{\mathcal{N}(Q)}\}$.
Consequently, we have
\begin{equation*}\label{equ:the alg01}
  C^*\{Q,I\}=C^*\{P_{\mathcal{R}(Q)},P_{\mathcal{N}(Q)}\}.
\end{equation*}
For an alternative proof of this equality, see \cite[Theorem~5.1]{BS}.
\end{remark}

Next, we recall  a new type of block matrix representation for a non-projection idempotent,  established in  \cite{TXF04}, which plays a crucial role in our analysis of the structure of the $C^*$-algebra generated by such an idempotent.

We first introduce some necessary notation.

\begin{definition}\label{defn of 6 operators and modules}  Let $Q\in \mathbb{B}(H)$ be an idempotent. Define
\begin{align}\label{eqn:defn of H1 and H2}&H_1=\mathcal{R}[m(Q)]\cap\mathcal{R}(Q), \quad H_2=\mathcal{R}[m(Q)]\cap\mathcal{N}(Q),\\
 \label{eqn:defn of H3 and H4}& H_3=\mathcal{N}[m(Q)]\cap\mathcal{R}(Q),  \quad H_4=\mathcal{N}[m(Q)]\cap\mathcal{N}(Q) \\
 \label{eqn:defn of H5 and H6}&H_5=\overline{\mathcal{R}\big[m(Q)Q\big(I-m(Q)\big)\big]}, \quad  H_6=\overline{\mathcal{R}\big[\big(I-m(Q)\big)Qm(Q)\big]}.
\end{align}
\end{definition}

A characterization of these subspaces of $H$ is given as follows.
\begin{proposition}
Let $Q\in \mathbb{B}(H)$ be an idempotent. Then
\begin{align}\label{equ:H1-2-H4} & H_1=\mathcal{R}(Q)\cap \mathcal{R}(Q^*),\quad  H_2= H_3=\{0\}, \quad H_4=\mathcal{N}(Q)\cap \mathcal{N}(Q^*),
\end{align}
where
$H_i(1\leq i\leq4)$ are defined in \eqref{eqn:defn of H1 and H2} and \eqref{eqn:defn of H3 and H4}.
\end{proposition}
\begin{proof} By \cite[Corollary~3.15]{TXF02}, we have $ H_2= H_3=\{0\}$.  Moreover, it is shown in  \cite[Remark~3.10]{TXF02} that
\begin{align}\label{equ:range4H1&H4}\mathcal{R}\big[m(Q)\big]=\mathcal{R}(|Q^*|+Q^*), \quad  m(Q)Q=\frac{1}{2}(|Q|+Q).
\end{align}

Let $x\in  \mathcal{R}(Q^*)\cap\mathcal{R}(Q)$. Then
$Q^*x=Qx=x$, and hence
\begin{align*}(I+|Q^*|)(I-|Q^*|)x=(I-QQ^*)x=0.\end{align*}
Since  $I+|Q^*|$ is invertible, we have  $(I-|Q^*|)x=0$.  Therefore,
$ x=Q^*x=|Q^*|x$,
which implies $x=\frac12(|Q^*|+Q^*)x$.
From the first equation in \eqref{equ:range4H1&H4}, we conclude that $x\in \mathcal{R}[m(Q)]$. Since  $x$ was arbitrary, it follows that
$$\mathcal{R}(Q^*)\cap\mathcal{R}(Q)\subseteq \mathcal{R}[m(Q)]\cap\mathcal{R}(Q)=H_1.$$

Conversely, let $x\in H_1$. By the second identity in \eqref{equ:range4H1&H4},
$$x=m(Q)Qx=\frac12(|Q|+Q)x=\frac12|Q|x+\frac12x,$$
which implies
$x=|Q|x\in \mathcal{R}(Q^*)$.
Since $x\in \mathcal{R}(Q)$ by assumption, it follows that
$x\in \mathcal{R}(Q^*)\cap\mathcal{R}(Q)$.  By the arbitrariness of $x$, we conclude that  $H_1\subseteq \mathcal{R}(Q^*)\cap\mathcal{R}(Q)$.
This confirms the first identity in \eqref{equ:H1-2-H4}.

By Lemma~\ref{lem:two relationships with m q}, we have
$ m(I-Q)=I-m(Q)$.
Therefore,  the last identity in \eqref{equ:H1-2-H4} follows by replacing $Q$ with $I-Q$.
\end{proof}

Given an  idempotent $Q\in \mathbb{B}(H)$, let $H_i(1\leq i\leq 6)$ be defined by \eqref{eqn:defn of H1 and H2}--\eqref{eqn:defn of H5 and H6}. Since every closed linear subspace  of a Hilbert space is always orthogonally complemented, the matched pair $\big(m(Q),Q\big)$ is always harmonious in the sense of \cite[Definition~4.1]{TXF05}.
By \cite[Lemma~4.1 and Remark~4.2]{TXF05},  $H_i\bot H_j$ for $i\ne j$, and
\begin{equation*}
  m(Q)=P_{H_1}+P_{H_2}+P_{H_5},\quad  I-m(Q)=P_{H_3}+P_{H_4}+P_{H_6}.
\end{equation*}
In view of  \eqref{equ:H1-2-H4}, we have $P_{H_2}=P_{H_3}=0$, so these expressions simplify to
\begin{equation*}\label{equ:P=p1+p2+p5}
  m(Q)=P_{H_1}+P_{H_5},\quad  I-m(Q)=P_{H_4}+P_{H_6}.
\end{equation*}
Let $U_{m(Q),Q}:H\rightarrow H_1\oplus H_4\oplus H_5\oplus H_6$ be  the unitary operator defined by
\begin{equation*}\label{equ:defn of Upq}U_{m(Q),Q}(x)=\Big(P_{H_1}(x),P_{H_4}(x),P_{H_5}(x),P_{H_6}(x)\Big)^T, \quad x\in H. \end{equation*}

We now state a key lemma, giving a block matrix representation for a matched pair $\big(m(Q),Q\big)$.
\begin{lemma}\label{lem:TXF}{\rm\cite[Theorem~4.3]{TXF04}}\label{lem:new representation for a general idempotent} Suppose that $Q\in\mathbb{B}(H)$ is a non-projection idempotent. Let  $H_1$, $H_4$, $H_5$ and $H_6$ be defined by \eqref{eqn:defn of H1 and H2}--\eqref{eqn:defn of H5 and H6}. Then
\begin{align*}
   & U_{m(Q),Q}\cdot m(Q)\cdot U_{m(Q), Q}^*=I_{H_1}\oplus 0_{H_4}\oplus I_{H_5}\oplus 0_{H_6}, \\
    & U_{m(Q),Q}\cdot Q\cdot U_{m(Q),Q}^*=I_{H_1}\oplus 0_{H_4}\oplus\widehat{Q},
\end{align*}
where
\begin{equation*}\label{eq:decomposition of Q0 hat2}
\widehat{Q}=\left(
                        \begin{array}{cc}
                          D & -\ell(D)U^* \\
                          U\ell(D) & U(I_{H_5}-D)U^* \\
                        \end{array}
                      \right)\in \mathbb{B}(H_5\oplus H_6),
\end{equation*}
in which $U\in \mathbb{B}(H_5, H_6)$ is a unitary operator,  $D\in \mathbb{B}(H_5)$ satisfies
\begin{equation*}\label{equ:the condition of A}
 D\geq I_{H_5},\quad \overline{\mathcal{R}(D^2-D)}=H_5,
\end{equation*}
and the function $\ell(t)$ is defined in \eqref{defn of lt}.
\end{lemma}

\begin{remark}
 The above condition on the operator $D$ can be rewritten as
\begin{equation}\label{equ:2propofD}
 D\geq I_{H_5},\quad \mathcal{N}(D-I_{H_5})=\{0\}.
\end{equation}
Moreover, since $D\ge I_{H_5}$, we have
\begin{equation}\label{equ:newtemp-003}\ell(D)=D^{\frac12}(D-I_{H_5})^{\frac12},\quad  \|\ell(D)\|=\sqrt{\|D\|(\|D\|-1)}=\ell(\|D\|).
\end{equation}
\end{remark}

\begin{remark}\label{rem:Wus point}If we shift the focus to the unitary equivalence of operators rather than the explicit construction of the unitary operators themselves,  an alternative proof of Lemma~\ref{lem:TXF} can be given as follows.

 Let $Q\in \mathbb{B}(H)$ be a non-projection idempotent. Since $Q$ is a quadratic operator with $\sigma(Q)=\{0,1\}$,  by \cite[Theorem~1.1]{TW},
   there exist Hilbert spaces $K_1,K_2,K_3$ and a unitary operator
$U_1: H\to K_1\oplus K_2\oplus K_3\oplus K_3$ such that
\begin{equation*}\label{block for idempotent}U_1QU_1^*=I_{K_1}\oplus 0_{K_2}\oplus \left(
                                                    \begin{array}{cc}
                                                      I_{K_3} & C \\
                                                      0_{K_3} & 0_{K_3} \\
                                                    \end{array}
                                                  \right),\end{equation*}
where $C\in\mathbb{B}(K_3)$ is a positive operator with $\mathcal{N}(C)=\{0\}$. Let
$$B=(I_{K_3}+C^2)^\frac12,\quad D=\frac12(I_{K_3}+B).$$
Then $D\ge I_{K_3}$ and
\begin{equation}\label{2exp of C}2\ell(D)=(B^2-I_{K_3})^\frac12=C.\end{equation} A direct computation yields
$$D-I_{K_3}=\big(2(B+I_{K_3})\big)^{-1}C^2,$$
which implies $$\mathcal{N}(D-I_{K_3})=\mathcal{N}(C^2)=\mathcal{N}(C)=\{0\}.$$
Since $U_1m(Q)U_1^*=m(U_1QU_1^*)$ and the operator $C$ can be expressed as \eqref{2exp of C},
by \cite[Theorems~3.1 and 3.4]{TXF02}, we have
\begin{equation*}
  U_1m(Q)U_1^*=I_{K_1}\oplus 0_{K_2}\oplus T,
\end{equation*}
where
\begin{align*}T=&\frac12\left(\begin{array}{cc}
        (B+I_{K_3})B^{-1} &  B^{-1}C\\
        CB^{-1} &  C\big[B(B+I_{K_3})\big]^{-1}C \\
      \end{array}\right)\\
=&\frac12\left(\begin{array}{cc}
        (B+I_{K_3})B^{-1} &  2B^{-1}\ell(D)\\
         2B^{-1}\ell(D) &  B^{-1}(B-I_{K_3})\\
      \end{array}\right)\\
=&\left(
    \begin{array}{cc}
      B^{-1}D & B^{-1}\ell(D) \\
      B^{-1}\ell(D) & B^{-1}(D-I_{K_3}) \\
    \end{array}
  \right).
\end{align*}

Now, define $W\in\mathbb{B}(K_3\oplus K_3)$ by
\begin{equation}\label{equ:unitary symmetry W}W=\left( \begin{array}{cc} B^{-\frac{1}{2}}D^{\frac12} & B^{-\frac{1}{2}}(D-I_{K_3})^{\frac12} \\ B^{-\frac{1}{2}}(D-I_{K_3})^{\frac12} & -B^{-\frac{1}{2}}D^{\frac12} \\  \end{array}  \right).\end{equation}
Then $W=W^*$ and $W^2=I_{K_3\oplus K_3}$, so $W$ is a unitary symmetry. A direct computation shows that $WTW=I_{K_3}\oplus 0_{K_3}$, and
\begin{equation}\label{equ:W4norm computation}W\left(
     \begin{array}{cc}
       I_{K_3} & 2\ell(D) \\
       0 & 0 \\
     \end{array}
   \right)W=\left(
              \begin{array}{cc}
                D & -\ell(D) \\
                \ell(D) & I_{K_3}-D \\
              \end{array}
            \right).\end{equation}

Finally, let $U_2=I_{K_1}\oplus I_{K_2}\oplus W$. Then $U_2\in \mathbb{B}(K_1\oplus K_2\oplus K_3\oplus K_3)$ is a unitary operator satisfying
\begin{align}\label{equ:mqu21}&U_2U_1m(Q)U_1^*U_2^*=I_{K_1}\oplus 0_{K_2}\oplus I_{K_3}\oplus 0_{K_3},\\
\label{equ:qu112}  &U_2U_1QU_1^*U_2^*=I_{K_1}\oplus 0_{K_2}\oplus \left(
              \begin{array}{cc}
                D & -\ell(D) \\
                \ell(D) & I_{K_3}-D \\
              \end{array}
            \right).
  \end{align}
Therefore, up to the unitary equivalence induced by the unitary operator $U_2U_1$, we obtain the desired decompositions of $m(Q)$ and $Q$.
\end{remark}

We now return to Lemma~\ref{lem:TXF}. To simplify  notation, define a unitary operator $V_Q\in \mathbb{B}(H, H_1\oplus H_4\oplus H_5\oplus H_5)$ by
\begin{equation}\label{equ:the unitary of V}
V_Q=\left(I_{H_1}\oplus I_{H_4}\oplus I_{H_5}\oplus U^*\right)U_{m(Q), Q}.
\end{equation}
 Then $m(Q)$ and $Q$  can be expressed as
\begin{align}\label{equ:the unitarily of M(Q)}
    &m(Q)=V_Q^*\big(I_{H_1}\oplus 0_{H_4}\oplus I_{H_5}\oplus 0_{H_5}\big)V_Q,  \\
 \label{equ:the unitarily of Q}  &Q=V_Q^*\Big(I_{H_1}\oplus 0_{H_4}\oplus \widetilde{Q}\Big) V_Q,
\end{align}
where $V_Q$ is a unitary operator defined by \eqref{equ:the unitary of V}, and
\begin{equation}\label{eq:new decomposition of Q2 hat2}
\widetilde{Q}=\left(
                        \begin{array}{cc}
                          D & -\ell(D) \\
                          \ell(D) & I_{H_5}-D \\
                        \end{array}
                      \right)\in \mathbb{B}(H_5\oplus H_5),
\end{equation}
in which $D\in \mathbb{B}(H_5)$ satisfies  \eqref{equ:2propofD}, and $\ell(D)$ is given by \eqref{equ:newtemp-003}.

Let $Q\in\mathbb{B}(H)$ be a non-projection idempotent. In the remainder of this subsection, we characterize all elements of the $\mathrm{C}^*$-algebra $C^*\{Q\}$ by representing them as $4 \times 4$ block matrices. This characterization is completed in Theorems~\ref{the:1 not in sigmaD} and \ref{th:1 in sigma}.

\begin{theorem}\label{the:1 not in sigmaD}
  Let  $Q\in \mathbb{B}(H)$ be a non-projection idempotent  admitted to the decomposition \eqref{equ:the unitarily of Q}, where $V_Q$ and $\widetilde{Q}$ are given by
  \eqref{equ:the unitary of V} and \eqref{eq:new decomposition of Q2 hat2}, respectively, such that the operator  $D$ in \eqref{eq:new decomposition of Q2 hat2} satisfies  \eqref{equ:2propofD}. If $1\notin \sigma(D)$, then
   \begin{equation*}\label{equ:Cstar of Q-2}
  C^*\{Q\}=\Big\{V_Q^*\big( \alpha I_{H_1}\oplus{0_{H_4}}\oplus T\big)V_Q: \alpha\in \mathbb{C}, T\in M_2\big(C^*\{D\}\big)\Big\}.
\end{equation*}
\end{theorem}
\begin{proof} Denote $0_{H_1}$, $0_{H_4}$ and $0_{H_5}$ simply by $0$. Let $d=\|D\|$, and let $\Gamma$ denote the set of polynomials vanishing at 0.

We aim to prove $ C^*\{V_QQV_Q^*\}=\mathfrak{B}$, where $\mathbb{C} I_{H_1}=\{\alpha I_{H_1}:\alpha\in \mathbb{C}\}$, and $\mathfrak{B}$ is the $C^*$-subalgebra of $\mathbb{B}(H_1\oplus H_4\oplus H_5\oplus H_5)$ given by
    \begin{equation*}\label{equ:Cstar of Q-22}\mathfrak{B}= \mathbb{C} I_{H_1}\oplus{0}\oplus  M_2\big(C^*\{D\}\big).
\end{equation*}
It is immediate from \eqref{equ:the unitarily of M(Q)} and \eqref{equ:the unitarily of Q} that
\begin{align}\label{Vm(q)V}
   & V_Qm(Q)V_Q^*=I_{H_1}\oplus 0\oplus I_{H_5}\oplus 0, \\
   \label{VqV}& V_QQV_Q^*=I_{H_1}\oplus 0\oplus \widetilde{Q}.
\end{align}  These identities imply
$C^*\{V_QQV_Q^*\}\subseteq \mathfrak{B}$, and
\begin{equation*}V_Qm(Q)Qm(Q)V_Q^*=I_{H_1}\oplus 0\oplus D\oplus 0.\end{equation*}
Therefore,
\begin{equation}\label{q mq q}p\big(V_Qm(Q)Qm(Q)V_Q^*\big)=p(1)I_{H_1} \oplus{0}\oplus  p(D)\oplus 0,\quad \forall p\in \Gamma,\end{equation}
which implies that
 \begin{equation}\label{equ:pd}p(1)I_{H_1} \oplus{0}\oplus  p(D)\oplus 0\in C^*\{V_QQV_Q^*\},\quad\forall p\in\Gamma,
\end{equation}
since  by  Theorem~\ref{thm:new express of Q} we have  $m(Q)\in C^*\{Q\}$.

Now, we establish the reverse inclusion $\mathfrak{B}\subseteq C^*\{V_QQV_Q^*\}$.
From the invertibility of $D$, we obtain $D^\dag=D^{-1}$. So by Lemma~\ref{lem: A dag contained sub alg}, $I_{H_5}=DD^\dag\in C^*\{D\}$. Thus, $C^*\{D\}=C^*\{D,I_{H_5}\}$. Hence, $C^*\{D\}$  is $C^*$-algebra isomorphic to $C\big(\sigma(D)\big)$ via the Gelfand transform.

For any $\alpha \in \mathbb{C}$ and $D_{11}\in C^*\{D\}$, let $A\in \mathfrak{B}$ be as
\begin{equation*}
 A=\alpha I_{H_1} \oplus{0}\oplus  D_{11} \oplus  0.
\end{equation*}
Write $D_{11}=f(D)$ for some $f\in C\big(\sigma(D)\big)$. As $D\ge I_{H_5}$ and $1\notin \sigma(D)$, we have $\sigma(D)\subseteq [1+\delta,d]$ for some $\delta>0$.
Define $\Omega=\{0,1\}\cup \sigma(D)$, and let $g$ be the unique extension of $f$ to  $\Omega$ satisfying $g(0)=0$ and $g(1)=\alpha$. Then $g\in C(\Omega)$, and
since $\Omega$ is a closed subset of the compact Hausdorff space $[0,d]$,  $g$ can be further extended to a function $h\in C[0,d]$. Thus,
$$A=h(1)I_{H_1} \oplus{0}\oplus  h(D)\oplus  0, \quad h\in C[0,d],\quad h(0)=0.$$
Take a sequence $\{p_n\}_{n=1}^{\infty}\subseteq \Gamma$ converging uniformly to $h$ on $[0,d]$. By \eqref{equ:pd}, each $A_n=p_n(1)I_{H_1} \oplus{0}\oplus  p_n(D)\oplus  0$ lies in $C^*\{V_QQV_Q^*\}$. Since $A_n\to A$ as $n\to\infty$, we conclude  $A\in C^*\{V_QQV_Q^*\}$. This shows that
\begin{equation}\label{eq:1general}\alpha I_{H_1} \oplus{0}\oplus  D_{11} \oplus  0\in C^*\{V_QQV_Q^*\},\quad  \forall \alpha\in \mathbb{C}, D_{11}\in C^*\{D\}.
 \end{equation}

Next, we prove that
\begin{equation}\label{eq:B12 in C*(Q)}
0 \oplus0\oplus  \begin{pmatrix}
  0 & D_{12} \\
  0 & 0
\end{pmatrix}\in C^*\{V_QQV_Q^*\},\quad \forall D_{12}\in C^*\{D\}.
\end{equation}
To this end, define $$D_{12}^{'}=-D_{12}D^{-\frac12}(D-I_{H_5})^{-\frac12}.$$  As shown earlier, we may write $D_{12}^{'}=\varphi(D)$ for some $\varphi\in C[0,d]$ satisfying $\varphi(0)=0$. Furthermore, there exists a sequence $\{q_n\}_{n=1}^{\infty}$ in $\Gamma$ such that $\|q_n(D)-D_{12}^{'}\|\to 0$ as $n\to \infty$.
For each $n\in \mathbb{N}$, define \begin{equation*}\label{equ:Tn}
      B_n=q_n\Big(V_Qm(Q)Qm(Q)V_Q^*\Big)\cdot V_Qm(Q)Q\big(I-m(Q)\big)V_Q^*.
    \end{equation*}
Using  \eqref{Vm(q)V}, \eqref{VqV} and Theorem~\ref{thm:new express of Q}, we obtain
\begin{equation}\label{equ:m(q)q(I-m(Q)}0\oplus 0\oplus \left(
                                                \begin{array}{cc}
                                                  0 & -\ell(D) \\
                                                  0 & 0 \\
                                                \end{array}
                                              \right)=
  V_Qm(Q)Q\big(I-m(Q)\big)V_Q^*\in C^*\{V_QQV_Q^*\}.
\end{equation}
Combining this with the definition of $B_n$ and Equation \eqref{q mq q}, we compute:
\begin{align*}
    B_n= 0\oplus 0\oplus \left(
                                                \begin{array}{cc}
                                                  0 & -q_n(D)\ell(D) \\
                                                  0 & 0 \\
                                                \end{array}
                                              \right)  \to 0\oplus 0\oplus \left(
                                                \begin{array}{cc}
                                                  0 & D_{12} \\
                                                  0 & 0 \\
                                                \end{array}
                                              \right)\quad \text{as}\ n \to \infty.
\end{align*}
Since $B_n\in C^*\{V_QQV_Q^*\}$ for all $n\in \mathbb{N}$, this establishes the validity of \eqref{eq:B12 in C*(Q)}.

For every $D_{21}\in C^*\{D\}$, a direct application of \eqref{eq:B12 in C*(Q)} yields
\begin{equation}\label{eq:B21 in C*(Q)}
0 \oplus{0}\oplus  \begin{pmatrix}
  0 & 0 \\
 D_{21} & 0
\end{pmatrix}=\left[0 \oplus{0}\oplus  \begin{pmatrix}
  0 & D_{21}^* \\
 0 & 0
\end{pmatrix}\right]^* \in C^*\{V_QQV_Q^*\}.
\end{equation}
For any $D_{22}\in C^*\{D\}$, substitute $D_{21}=I_{H_5}$ into \eqref{eq:B21 in C*(Q)} and $D_{12}=D_{22}$ into $\eqref{eq:B12 in C*(Q)}$; multiplying  the resulting expressions gives
\begin{equation}\label{eq:B22 in C*(Q)}
0 \oplus{0}\oplus  \begin{pmatrix}
  0 & 0 \\
  0 & D_{22}
\end{pmatrix} \in C^*\{V_QQV_Q^*\}.
\end{equation}
Therefore, the reverse inclusion  follows  from \eqref{eq:1general}, \eqref{eq:B12 in C*(Q)}, \eqref{eq:B21 in C*(Q)} and \eqref{eq:B22 in C*(Q)}.
This completes the proof.
\end{proof}

\begin{remark}Under the condition of Theorem~\ref{the:1 not in sigmaD}, let
$$T=\left(
      \begin{array}{cc}
        T_{11} & T_{12} \\
        T_{21}& T_{22} \\
      \end{array}
    \right)
\in M_2\big(C^*\{D\}\big).$$
For $1\le i,j\le 2$, write $T_{ij}=g_{ij}(D)$, where each $g_{ij}\in C\big(\sigma(D)\big)$. Given any complex number $\alpha$, since $\sigma(D)$ is a compact subset of $(1,\|D\|]$, there exists $f_{11}\in C[1,\|D\|]$ extending $g_{11}$ and satisfying $f_{11}(1)=\alpha$. For each $i,j$ with $i+j>2$, let $f_{ij}\in C[1,\|D\|]$ be an extension of $g_{ij}$ satisfying $f_{ij}(1)=0$. Then an alternative expression for $\alpha I_{H_1}\oplus 0_{H_4}\oplus T$ is
$$f_{11}(1)I_{H_1}\oplus 0_{H_4}\oplus \left(
                                   \begin{array}{cc}
                                   f_{11}(D) & f_{12}(D) \\
                                     f_{21}(D) & f_{22}(D) \\
                                   \end{array}
                                 \right).$$
\end{remark}

To facilitate the subsequent discussion, we  present the following definition.
\begin{definition}\label{def:Jd}
Suppose $D\in\mathbb{B}(K)_{\mbox{sa}}$ with $1\in\sigma(D)$. Let
$I_{\{1\}}$ be the closed  ideal of  $C\big(\sigma(D)\big)$ defined by
$$I_{\{1\}}=\big\{g\in C\big(\sigma(D)\big):g(1)=0\big\},$$
which induces a $C^*$-subalgebra of $M_2\Big(C\big(\sigma(D)\big)\Big)$ as
\begin{equation}
  \label{defn jd}J_D=\left\{(f_{ij})_{1\leq i,j\leq 2}\in M_2\Big(C\big(\sigma(D)\big)\Big):f_{ij}\in I_{\{1\}}\ \text{when}\ i+j>2 \right\}.
\end{equation}
Furthermore, for any $f=(f_{ij})_{1\leq i,j\leq 2}\in M_2\Big(C\big(\sigma(D)\big)\Big)$ and $t\in\sigma(D)$, let $f_t\in M_2(\mathbb{C})$ and $f(D)\in M_2\big(C^*\{D,I_K\}\big)$ be defined by
\begin{equation}\label{2ff}f_t=\left(
         \begin{array}{cc}
           f_{11}(t) & f_{12}(t) \\
           f_{21}(t) & f_{22}(t) \\
         \end{array}
       \right),\quad f(D)=\left(
         \begin{array}{cc}
           f_{11}(D) & f_{12}(D) \\
           f_{21}(D) & f_{22}(D) \\
         \end{array}
       \right).
\end{equation}
\end{definition}

\begin{theorem}\label{th:1 in sigma}
  Let  $Q\in \mathbb{B}(H)$ be a non-projection idempotent  admitted to the decomposition \eqref{equ:the unitarily of Q}, where $V_Q$ and $\widetilde{Q}$ are given by
  \eqref{equ:the unitary of V} and \eqref{eq:new decomposition of Q2 hat2}, respectively, such that the operator  $D$ in \eqref{eq:new decomposition of Q2 hat2} satisfies  \eqref{equ:2propofD}. If $1\in \sigma(D)$, then
   \begin{equation*}\label{equ:CstarQ-1in}
  C^*\{Q\}=\Big\{V_Q^*\Big(f_{11}(1) I_{H_1} \oplus{0_{H_4}}\oplus f(D)\Big)V_Q: f=(f_{ij})_{1\leq i,j\leq 2}\in J_D\Big\},
\end{equation*}
where $J_D$ is defined by \eqref{defn jd}.
\end{theorem}
\begin{proof}We aim to prove
$C^*\{V_QQV_Q^*\} = \mathfrak{A}$,
where
\begin{align*}
   & \mathfrak{A}=\big\{f_{11}(1) I_{H_1} \oplus{0}\oplus f(D):f=(f_{ij})_{1\leq i,j\leq 2}\in J_D\big\}.
\end{align*}
From the proof of Theorem~\ref{the:1 not in sigmaD}, we see  that $g(D)\in C^*\{D\}$ for any $g\in C\big(\sigma(D)\big)$.
Since $J_D$ is evidently a $C^*$-subalgebra of $M_2\Big(C\big(\sigma(D)\big)\Big)$, it follows that
$\mathfrak{A}$ is a $C^*$-subalgebra of $\mathbb{C}I_{H_1}\oplus{0}\oplus  M_2\big(C^*\{D\}\big)$.
By \eqref{VqV} and \eqref{eq:new decomposition of Q2 hat2}, $V_QQV_Q^*\in \mathfrak{A}$, which implies  $C^*\{V_QQV_Q^*\}\subseteq  \mathfrak{A}$.

We now establish the reverse inclusion.  Combining  \eqref{Vm(q)V} and \eqref{equ:pd}, we conclude that
$p(1)I_{H_1} \oplus{0}\oplus  p(D) \oplus 0 \in C^*\{V_QQV_Q^*\}$ for any polynomial $p$. This means
\begin{equation}\label{eq:f11 in C*(Q)}
 f_{11}(1)I_{H_1} \oplus{0}\oplus  f_{11}(D) \oplus 0 \in C^*\{V_QQV_Q^*\},\quad \forall f_{11}\in C\big(\sigma(D)\big).
\end{equation}
In particular,
\begin{equation}\label{equ:g}
  0\oplus 0\oplus -g(D)\ell(D) \oplus 0\in C^*\{V_QQV_Q^*\},\quad \forall g\in C\big(\sigma(D)\big).
\end{equation}

Given an arbitrary element $f_{12}$ in $I_{\{1\}}$, define $f_{12}'\in I_{\{1\}}$  by
$$f_{12}'(x)=\frac{f_{12}(x)}{x},\quad x\in\sigma(D).$$
Select a sequence  of polynomials $\{p_n\}_{n=1}^{\infty}$ such that $p_n(1)=0$ for all $n\in \mathbb{N}$, and $$\lim\limits_{n\to \infty}\max\limits_{x\in \sigma(D)}\big|p_n(x)-f_{12}'(x)\big|=0.$$
For each $n\in \mathbb{N}$, write $p_n(x)=(x-1)q_n(x)$, where $q_n(x)$ is also a  polynomial.  Since $\sigma(D)$ is  bounded in $\mathbb{C}$, it follows that
$$\lim\limits_{n\to \infty}\max\limits_{x\in \sigma(D)}|x(x-1)q_n(x)-f_{12}(x)|=\lim\limits_{n\to \infty}\max\limits_{x\in \sigma(D)}|x\big(p_n(x)-f_{12}'(x)\big)|=0,$$
and hence \begin{equation*}\label{equ:norm of D}
                  \lim\limits_{n\to \infty}\|\ell^2(D)q_n(D)-f_{12}(D)\|=0.
                \end{equation*}
By \eqref{equ:m(q)q(I-m(Q)}, this yields
$$\lim\limits_{n\to \infty}S_nV_Qm(Q)Q\big(I-m(Q)\big)V_Q^*=0\oplus 0\oplus \left(
                                                \begin{array}{cc}
                                                  0 & f_{12}(D) \\
                                                  0 & 0 \\
                                                \end{array}
                                              \right), $$
where  $S_n= 0\oplus 0\oplus -q_n(D)\ell(D) \oplus 0$ for each $n\in\mathbb{N}$. Combining this limit  with \eqref{equ:g} and \eqref{equ:m(q)q(I-m(Q)},
we conclude that
\begin{equation}\label{equ:f12 in}
  0\oplus 0\oplus \left(
                                                \begin{array}{cc}
                                                  0 & f_{12}(D) \\
                                                  0 & 0 \\
                                                \end{array}
                                              \right)\in C^*\{V_QQV_Q^*\},\quad \forall f_{12}\in I_{\{1\}}.
\end{equation}
A similar argument, analogous to the proof of Theorem~\ref{the:1 not in sigmaD}, shows that
\begin{equation}\label{equ:f21 in}
  0\oplus 0\oplus \left(
                                                \begin{array}{cc}
                                                  0 & 0 \\
                                                  f_{21}(D) & 0 \\
                                                \end{array}
                                              \right)\in C^*\{V_QQV_Q^*\}, \quad \forall f_{21}\in I_{\{1\}}.
\end{equation}

Let $\mathfrak{B}=\big\{h(D):h\in I_{\{1\}}\big\}$.
Since $I_{\{1\}}$ is a closed two-sided ideal of $C\big(\sigma(D)\big)$, $\mathfrak{B}$ is a $C^*$-subalgebra of $C^*\{D\}$. Let $\{e_\lambda\}_{\lambda\in\Lambda}$ be an arbitrary approximate identity for  $\mathfrak{B}$. For any $f_{22}\in I_{\{1\}}$,
$$0\oplus 0\oplus 0\oplus f_{22}(D)=\lim\limits_{\lambda\in \Lambda}  \left[0\oplus 0\oplus \left(
                                                \begin{array}{cc}
                                                  0 & 0 \\
                                                  e_\lambda & 0 \\
                                                \end{array}
                                              \right)\right]\left[0\oplus 0\oplus \left(
                                                \begin{array}{cc}
                                                  0 & f_{22}(D) \\
                                                  0 & 0 \\
                                                \end{array}
                                              \right)\right].$$
By \eqref{equ:f21 in} and \eqref{equ:f12 in}, we conclude that
\begin{equation}\label{equ:f22 in}
  0\oplus 0\oplus 0\oplus f_{22}(D)\in C^*\{V_QQV_Q^*\}.
\end{equation}
The reverse inclusion now follows  from \eqref{eq:f11 in C*(Q)} and  \eqref{equ:f12 in}--\eqref{equ:f22 in}.
\end{proof}

\begin{remark}Let $Q\in\mathbb{B}(H)$ be a non-projection idempotent. Combining the last equality in \eqref{equ:H1-2-H4} with Theorems~\ref{the:1 not in sigmaD} and \ref{th:1 in sigma}, we find that
$$I\in C^*\{Q\}\Longleftrightarrow \mathcal{N}(Q)\cap \mathcal{N}(Q^*)=\{0\}\ \text{and}\ 1\notin \sigma(D).$$
Also, since $C^*\{Q,I\}=\big\{\lambda I+x:\lambda \in\mathbb{C}, x\in C^*\{Q\}\big\}$, the $4\times 4$ block matrix representations for elements in  $C^*\{Q,I\}$ can be derived from Theorems~\ref{the:1 not in sigmaD} and \ref{th:1 in sigma}.

\end{remark}

\subsection{The universal $r$-idempotents}\label{sec:universal-r-idempotents}

It is known that for any idempotent $Q$ on a Hilbert space, $Q$ is a projection if and only if $\|Q\|\le 1$. Since we are concerned with non-projection idempotents, we only consider the case where $r>1$ in the definition below.

\begin{definition}\label{defn:universalr}Let $r>1$. A non-projection idempotent $Q$ in a $C^*$-algebra is called a universal $r$-idempotent if  $\|Q\|=r$,  and it satisfies the following universal property:

For every  non-projection idempotent $Q'$ in any  $C^*$-algebra with
$\|Q'\| \leq r$, there exists a $C^*$-algebra homomorphism $\pi: C^*\{Q\} \to C^*\{Q'\}$ such that  $\pi(Q) = Q'$.
\end{definition}

\begin{remark} By definition, if  $Q_1$ and $Q_2$ are both  universal $r$-idempotents, the $C^*$-algebras they generate are canonically isomorphic.
Since every $C^*$-algebra can be embedded as a $C^*$-subalgebra of $\mathbb{B}(K)$ for some Hilbert space $K$, we may restrict our attention to idempotents acting on Hilbert spaces.
\end{remark}

To offer  an alternative characterization of the universal $r$-idempotents, we present the following lemma.

\begin{lemma}\label{lem:4universalr01}Let $Q_1\in\mathbb{B}(H)$ and $Q_2\in\mathbb{B}(K)$ be idempotents, and let $\pi:C^*\{Q_1\}\to C^*\{Q_2\}$ be a $C^*$-algebra homomorphism. Then the following statements are equivalent:
\begin{enumerate}
  \item [{\rm (i)}] $\pi(Q_1)=Q_2$;
  \item [{\rm (ii)}] $\pi\big(m(Q_1)\big)=m(Q_2)$ and $\pi(P_{\mathcal{R}(Q_1)})=P_{\mathcal{R}(Q_2)}$.
\end{enumerate}
\end{lemma}
\begin{proof}We establish the equivalence by proving both implications: (i)$\Longrightarrow$(ii) and (ii)$\Longrightarrow$(i).

(i)$\Longrightarrow$(ii). Assume $\pi(Q_1)=Q_2$. By \cite[Theorem~3.6]{TXF02},
$\pi\big(m(Q_1)\big)=m\big(\pi(Q_1)\big)=m(Q_2)$. It remains to prove $\pi(P_{\mathcal{R}(Q_1)})=P_{\mathcal{R}(Q_2)}$.
We split the proof into two cases:

\textbf{Case 1:}\ $I_H\in C^*\{Q_1\}$. It is easily seen that
$\pi\big(C^*\{Q_1\}\big)=C^*\{Q_2\}$, so $\pi(I_H)$ is the unit of $C^*\{Q_2\}$. Restricting $K$ to the closed subspace $\pi(I_H)K$, we may assume  $\pi:C^*\{Q_1\}\to C^*\{Q_2\}$ is unital. Applying $\pi$ to both sides of \eqref{equ:projection Prq} gives $\pi(P_{\mathcal{R}(Q_1)})=P_{\mathcal{R}(Q_2)}$.

\textbf{Case 2:}\ $I_H\notin C^*\{Q_1\}$. Extend $\pi$ to  a unital $C^*$-algebra homomorphism $\widetilde{\pi}: C^*\{Q_1,I_H\}\to C^*\{Q_2,I_K\}$ by defining
$$\widetilde{\pi}(\lambda I_H+x)=\lambda I_K+\pi(x),\quad \lambda\in \mathbb{C}, x\in C^*\{Q_1\}.$$ Applying the result from Case 1 to $\widetilde{\pi}$, we obtain
\begin{equation*}\pi\big(P_{\mathcal{R}(Q_1)}\big)=\widetilde{\pi}\big(P_{\mathcal{R}(Q_1)}\big)=P_{\mathcal{R}\big(\widetilde{\pi}(Q_1)\big)}=P_{\mathcal{R}\big(\pi(Q_1)\big)}=P_{\mathcal{R}(Q_2)}.
\end{equation*}

(ii)$\Longrightarrow$(i). For $i=1,2$, decompose the operators as  follows:
\begin{align*}&P_{\mathcal{R}(Q_i)}\big(2m(Q_i)-I\big)P_{\mathcal{R}(Q_i)}=2P_{\mathcal{R}(Q_i)}m(Q_i)P_{\mathcal{R}(Q_i)}-P_{\mathcal{R}(Q_i)},\\
&P_{\mathcal{R}(Q_i)}\big(2m(Q_i)-I\big)=2P_{\mathcal{R}(Q_i)}m(Q_i)-P_{\mathcal{R}(Q_i)}.
\end{align*}
Using identity \eqref{equ:new express of Q} and Lemma~\ref{lem:keepMP}, we thus  obtain  $\pi(Q_1)=Q_2$.
\end{proof}

Utilizing Definition~\ref{defn:universalr} and Lemma~\ref{lem:4universalr01}, we provide an equivalent characterization of universal r-idempotents.

\begin{proposition}\label{prop:alternativeup}Let $Q$ be a non-projection idempotent in a $C^*$-algebra with $\|Q\|=r$. Then $Q$ is a universal $r$-idempotent if and only if for every  non-projection idempotent $Q'$ in any $C^*$-algebra satisfying $\|Q'\| \leq r$, there exists a $C^*$-algebra homomorphism $\pi: C^*\{Q\} \to C^*\{Q'\}$ such that
$$\pi\big(m(Q)\big)=m(Q')\quad\mbox{and}\quad \pi(P_{\mathcal{R}(Q)})=P_{\mathcal{R}(Q')}.$$
\end{proposition}

\begin{remark}\label{rem:differup}Let $Q$ be a universal $r$-idempotent. It is notable that  the above  universal property associated with the pair of projections $m(Q)$ and $P_{\mathcal{R}(Q)}$ is distinct from  the one described in \cite[Proposition~1.1]{RS}. As observed in \cite[Remark~3.3]{TXF02}, $m(Q)$ and $P_{\mathcal{R}(Q)}$ are homotopy equivalent for any idempotent $Q$. Consequently,  if $P_1,P_2\in\mathbb{B}(K)$ are projections that are not homotopy equivalent, then there exists no $C^*$-algebra homomorphism $\pi$ such that
$\pi\big(m(Q)\big)=P_1$ and $\pi(P_{\mathcal{R}(Q)})=P_2$.
\end{remark}

To construct the universal $r$-idempotent $Q_r$ below,  we require the following lemma.
\begin{lemma}\label{lem:useful norm01}For every $D\in \mathbb{B}(K)$ with $D\geq I_K$, define the idempotent  $\widetilde{Q}\in \mathbb{B}(K\oplus K)$ by replacing $H_5$ with $K$ in \eqref{eq:new decomposition of Q2 hat2}. Then  $\|\widetilde{Q}\|=2\|D\|-1$.
\end{lemma}
\begin{proof}Let $B=2D-I_K$, and define the unitary symmetry $W\in \mathbb{B}(K \oplus K)$ by replacing $K_3$ with $K$ in \eqref{equ:unitary symmetry W}. By \eqref{equ:W4norm computation}, we have
 \begin{equation*}\label{equ:4quadricform}
\widetilde{Q}=W\left( \begin{array}{cc}
    I_K & 2\ell(D) \\
    0 & 0 \\
  \end{array}
\right)W.
\end{equation*}
Therefore,
$$\|\widetilde{Q}\|=\sqrt{\|\widetilde{Q}\widetilde{Q}^*\|}=\|2D-I_{K}\|.$$
Since $D\ge I_K$, its spectrum $\sigma(D)$ is contained in $[1,\|D\|]$. Thus,  $$\|\widetilde{Q}\|=\sup\limits_{t\in\sigma(D)}(2t-1)=2\|D\|-1.$$
This completes the proof.
\end{proof}

\begin{definition}
For  $r>1$, let $d=\frac12(r+1)$.  Define the idempotent $Q_r$ in the $C^*$-algebra $\mathbb{C}\oplus \{0\}\oplus M_2\big(C[1,d]\big)$ by
\begin{equation}\label{equ:Qr2}
  Q_r=1\oplus 0\oplus\left(
                                                                  \begin{array}{cc}
                                                                    D_0 & -\ell(D_0) \\
                                                                   \ell(D_0) & I_0-D_0 \\
                                                                  \end{array}
                                                                \right),
\end{equation}
where $I_0,D_0$ and $\ell(D_0)$ are elements of $C[1,d]$ given by
\begin{equation}\label{equ:D0I0}I_0(t)\equiv 1,\quad D_0(t)=t, \quad \ell(D_0)(t)=\sqrt{t(t-1)},\quad  t\in [1,d].\end{equation}
\end{definition}

\begin{theorem}\label{th:C*-homorphism2}
For every $r> 1$, the idempotent $Q_r$ defined by \eqref{equ:Qr2} is a universal $r$-idempotent.
\end{theorem}
\begin{proof}Let $d=\frac12(r+1)$. As $r>1$, we have  $d>1$. We express  $Q_r$ as $Q_r=1\oplus 0\oplus \widetilde{Q_r}$, where
\begin{equation}\label{expression of Qrtilde}\widetilde{Q_r}=\left(
                        \begin{array}{cc}
                          D_0 & -\ell(D_0)\\
                          \ell(D_0) & I_{0}-D_0 \\
                        \end{array}
                      \right)\in M_2\big(C[1,d]\big).\end{equation}
By \eqref{equ:D0I0}, $D_0\ge I_0$, $\sigma(D_0)=[1,d]$ and  $\ell(D_0)=D_0^{\frac12}(D_0-I_{0})^{\frac12}$. Hence, $\|D_0\|=d$. Applying
Lemma~\ref{lem:useful norm01}, we find that
$$\|Q_r\|=\max\big\{1,2\|D_0\|-1\big\}=\max\left\{1,2d-1\right\}=2d-1=r.$$
Moreover, from the derivation of \eqref{equ:mqu21} and \eqref{equ:qu112} in Remark~\ref{rem:Wus point}, we deduce that $m(Q_r)=1\oplus 0\oplus I_0\oplus 0$. Since $1\in\sigma(D_0)$,
Theorem~\ref{th:1 in sigma}  implies
 \begin{equation}\label{equ:Cstar of Qr2}
  C^*\{Q_r\}   = \big\{f_{11}(1) \oplus{0}\oplus f:f=(f_{ij})_{1\leq i,j\leq 2}\in J_{D_0}\big\},
\end{equation}
where  $J_{D_0}$ is defined as in  \eqref{defn jd}  by replacing $D$ with $D_0$.

Given an arbitrary non-projection  idempotent $Q\in \mathbb{B}(H)$  with  $\|Q\|\le r$, represent $Q$ in the form given in \eqref{equ:the unitarily of Q}, where
$V_Q$ and $\widetilde{Q}$ are defined  by \eqref{equ:the unitary of V} and \eqref{eq:new decomposition of Q2 hat2}, respectively, such that the operator  $D$ in \eqref{eq:new decomposition of Q2 hat2} satisfies  \eqref{equ:2propofD}. It follows from Lemma~\ref{lem:useful norm01} that
\begin{equation*}\label{equ:the norm of D}
  \|D\|=\frac12(\|Q\|+1)\le\frac12(r+1)=d.
\end{equation*} Thus,
\begin{equation*}\label{inc:sigmaD}
  \sigma(D)\subseteq [1,\|D\|]=[1,d].
\end{equation*}
Hence, for every $g\in C[1,d]$, the functional calculus $g(D)$ is well-defined.

Now, let  $T$ be an element of  $C^*\{Q_r\}$, expressed in the form given by the right-hand side of \eqref{equ:Cstar of Qr2}. Define $\pi_Q(T)\in\mathbb{B}(H)$ by
\begin{equation}\label{equ:homo of Q}
  \pi_Q(T)=
V_Q^*\left[f_{11}(1)I_{H_1} \oplus{0}\oplus \big(f_{ij}(D)\big)_{1\le i,j\le 2}\right]V_Q.
\end{equation}
Regardless of whether $1\in \sigma(D)$, Theorems~\ref{the:1 not in sigmaD} and \ref{th:1 in sigma} imply that $\pi_Q(T)\in C^*\{Q\}$. One can easily verify that $\pi_Q:C^*\{Q_r\}\to  C^*\{Q\}$  is  a $C^*$-algebra homomorphism
satisfying  $\pi_Q(Q_r)=Q$.
\end{proof}

\begin{remark}\label{rem:4universal03} For $r>1$, let $\widetilde{Q_r}$ be the idempotent given by  \eqref{expression of Qrtilde}. From the proof of the preceding theorem, we see  that $\widetilde{Q_r}$ is also a universal $r$-idempotent.
\end{remark}

\begin{theorem}\label{thm:nscondition4ui} Let  $Q\in \mathbb{B}(H)$ be a non-projection idempotent  of the form given in \eqref{equ:the unitarily of Q} with $\|Q\|=r$, where $V_Q$ and $\widetilde{Q}$ are defined in
  \eqref{equ:the unitary of V} and \eqref{eq:new decomposition of Q2 hat2}, respectively, such that the operator  $D$ in \eqref{eq:new decomposition of Q2 hat2} satisfies  \eqref{equ:2propofD}. Then $Q$ is a universal $r$-idempotent if and only if $\sigma(D)=[1,\|D\|]$.
\end{theorem}
\begin{proof}Suppose $\sigma(D)=[1,\|D\|]$.  Then $1\in\sigma(D)$ and $\|D\|=\frac12(r+1)=d$. By Theorem~\ref{th:1 in sigma}, every element $\widetilde{T}\in C^*\{Q\}$ can be expressed as
$$\widetilde{T}=V_Q^*\Big(f_{11}(1) I_{H_1} \oplus{0_{H_4}}\oplus f(D)\Big)V_Q$$
for some $f=(f_{ij})_{1\leq i,j\leq 2}\in J_D$. By \eqref{equ:homo of Q}, we have
$\widetilde{T}=\pi_Q(T)$, where
$$T=f_{11}(1) \oplus 0 \oplus f(D_0)\in C^*\{Q_r\},$$
in which $Q_r$ is the idempotent defined in \eqref{equ:Qr2}. It follows that
$$\|f(D)\|=\max_{t\in [1,\|D\|]}\|f_t\|=\max_{t\in [1,d]}\|f_t\|=\|f(D_0)\|,$$
where $f_t$ is defined by \eqref{2ff}. Since
$\|f_1\|=\big\|\big(f_{ij}(1)\big)\big\|\ge |f_{11}(1)|$, we have
$$\|\pi_Q(T)\|=\max\{\|f_{11}(1)I_{H_1}\|, \|f(D)\|\}=\|f(D)\|=\|f(D_0)\|=\|T\|.$$
 Hence, the $C^*$-algebra homomorphism $\pi_Q:C^*\{Q_r\}\to  C^*\{Q\}$  is injective.
Notably, $\pi_Q$ is always surjective. Thus, $\pi_Q$ is a $C^*$-algebra isomorphism. Consequently, $Q$ is a universal $r$-idempotent, since $Q_r$ is.

Conversely, assume $Q$ is a universal $r$-idempotent. Then there exists a $C^*$-algebra homomorphism $\pi: C^*\{Q\}\to C^*\{Q_r\}$ such that $\pi(Q)=Q_r$.  Now suppose, for contradiction, that there exists $\lambda\in [1,d]$ such that $\lambda\notin\sigma(D)$. Define
$$T=\lambda m(Q)-m(Q)Qm(Q).$$
Then by Theorem~\ref{thm:new express of Q}, $T\in C^*\{Q\}$. Moreover, by \eqref{equ:the unitarily of M(Q)} and \eqref{equ:the unitarily of Q}, we have
$$T=V_Q^*\Big[(\lambda-1)I_{H_1}\oplus 0\oplus (\lambda I_{H_5}-D)\oplus 0\Big]V_Q,
$$
which implies that $T$ is M-P invertible. By Lemma~\ref{lem:keepMP},  $\pi(T)$ is also M-P invertible in $C^*\{Q_r\}$. From Lemma~\ref{lem:4universalr01}, we obtain
$$\pi(T)=(\lambda-1)\oplus 0\oplus (\lambda I_0-D_0)\oplus 0.$$
Since $D_0$ is self-adjoint,  $\lambda I_0-D_0$ is M-P invertible if and only if $\lambda\notin \sigma(D_0)$ or $\lambda$ is an isolated point of $\sigma(D_0)$. However,  $\sigma(D_0)=[1,d]$, so both conditions are violated. Hence, $\lambda I_0-D_0$ is not M-P invertible, contradicting the M-P invertibility of $\pi(T)$.
Therefore, $\sigma(D)=[1,d]$.
\end{proof}

Next, we  present a simplified characterization of universal $r$-idempotents.
\begin{theorem}\label{thm:4 r-idem nonblock}Let  $Q\in \mathbb{B}(H)$ be a non-projection idempotent with $\|Q\|=r$. Then $Q$ is a universal $r$-idempotent if and only if $\sigma(A)=[1,d]$,
where $d=\frac12(r+1)$, and the operator $A\in\mathbb{B}(H)$ is defined by
$$A=m(Q)Qm(Q)+I-m(Q).$$
\end{theorem}
\begin{proof}By \eqref{equ:the unitarily of M(Q)}--\eqref{eq:new decomposition of Q2 hat2} and Lemma~\ref{lem:useful norm01}, $\|D\|=d$, $\sigma(D)\subseteq [1,d]$, and
\begin{equation}\label{4 A r idem}A=V_Q^*\Big[I_{H_1}\oplus I_{H_4}\oplus D\oplus I_{H_5}\big]V_Q,\end{equation}
where $D\in\mathbb{B}(H_5)$  satisfies  \eqref{equ:2propofD}.  It follows from \eqref{4 A r idem}  that $A\ge I$, so
$$\sigma(A)\subseteq [1,\|A\|]=[1,\|D\|]=[1,d].$$
By \eqref{4 A r idem}, $\lambda I-A$ is invertible in $\mathbb{B}(H)$ if and only if $\lambda I_{H_5}-D$ is invertible in $\mathbb{B}(H_5)$ for any $\lambda\in (1,d]$. Since both $\sigma(D)$ and $\sigma(A)$ are closed subsets of $\mathbb{C}$, we conclude
$$[1,d]\subseteq \sigma(D)\Longleftrightarrow (1,d]\subseteq \sigma(D)\Longleftrightarrow   (1,d] \subseteq \sigma(A)\Longleftrightarrow   [1,d] \subseteq \sigma(A).$$
The desired result follows from Theorem~\ref{thm:nscondition4ui}.
\end{proof}

\begin{corollary}\label{cor:Q&I-Quri}Let  $Q\in \mathbb{B}(H)$ be a non-projection idempotent. Then $Q$ is a universal $r$-idempotent if and only if $I-Q$ is a universal $r$-idempotent.
\end{corollary}
\begin{proof}Suppose $Q$ is a universal $r$-idempotent. By Definition~\ref{defn:universalr}, $\|Q\|=r$, hence  $\|I-Q\|= r$, since $\|I-Q\|=\|Q\|$. By Theorem~\ref{thm:4 r-idem nonblock}, $I-Q$ is a universal $r$-idempotent if and only if $\sigma(B)=[1,d]$, where $d=\frac12(r+1)$, and $B\in\mathbb{B}(H)$ is given by
$$B=m(I-Q)(I-Q)m(I-Q)+I-m(I-Q).$$
Using Lemma~\ref{lem:two relationships with m q} and \eqref{equ:the unitarily of M(Q)}--\eqref{eq:new decomposition of Q2 hat2}, we compute:
\begin{align*}B=&\big(I-m(Q)\big)(I-Q)\big(I-m(Q)\big)+m(Q)\\
=&V_Q^*\Big[I_{H_1}\oplus I_{H_4}\oplus I_{H_5}\oplus D\big]V_Q.
\end{align*}
From this expression for $B$ together with \eqref{4 A r idem} and Theorem~\ref{thm:4 r-idem nonblock}, we conclude $\sigma(B)=\sigma(A)=[1,d]$. Thus, $I-Q$ is a universal $r$-idempotent.

Conversely, since $Q=I-(I-Q)$, the reverse implication follows immediately.
\end{proof}

\begin{remark}\label{rem:4identification}Let $H_0=L^2[1,d]$ with $d>1$.  For any $f\in C[1,d]$,  define the multiplication operator  $M_f\in \mathbb{B}(H_0)$ by
 \begin{equation*}\label{equ:Mf011}M_f(g)=fg,\quad g\in H_0.\end{equation*}
This induces  a faithful unital $C^*$-algebra homomorphism
   from $C[1,d]$ into $\mathbb{B}(H_0)$.
 For simplicity, throughout the rest of this paper, we identify each $f\in C[1,d]$ with its image $M_f$.
\end{remark}

\begin{theorem}\label{non-uniqueness universal}
For every $r>1$,  there exist two universal $r$-idempotents that are not unitarily equivalent.
\end{theorem}
\begin{proof}Let $d=\frac12(r+1)$, and let $\widetilde{Q_r}\in M_2\big(C[1,d]\big)$ be the universal $r$-idempotent defined by \eqref{expression of Qrtilde}.
Define an additional idempotent  $q_r$ by
\begin{equation*}\label{equ:Q3}
  q_r=\left(
                                                                  \begin{array}{cc}
                                                                    D_{1}& -\ell(D_{1}) \\
                                                                   \ell(D_{1}) & I_0-D_{1}\\
                                                                  \end{array}
                                                                \right)\in M_2\big(C[1,d]\big),
\end{equation*}
  where $I_0$ is defined by \eqref{equ:D0I0} and  $D_{1}:[1,d]\to\mathbb{C}$ is the piecewise function
  \begin{equation*}
    D_{1}(t)=\left\{
            \begin{array}{ll}
              1, & \hbox{if $t\in [1,\frac{d+1}{2}]$,} \\
              2t-d, & \hbox{if $t\in(\frac{d+1}{2},d]$.}
            \end{array}
          \right.
  \end{equation*}
Then $D_{1}\in C[1,d]$ and $\sigma(D_{1})=[1,d]$. From the proof of Theorem~\ref{thm:nscondition4ui},  it follows that $q_r$ is also a universal $r$-idempotent.

 Let $g\in H_0=L^2[1,d]$. By \eqref{equ:D0I0},   we have
 $[(D_0-I_0)g](t) =(t-1)g(t)$ for $t\in[1,d]$.
 Hence,  $g\in \mathcal{N}(D_0-I_0)\Longleftrightarrow g=0\ \mbox{a.e.}$, so $\mathcal{N}(D_0-I_0)=\{0\}$. However,
 this does not hold when $D_0$ is replaced by  $D_1$. Indeed, by definition of $D_1$, $[(D_1-I_0)h](t)=0$ for all $t\in [1,d]$, where $h$ is the characteristic function of  the subset $[1,\frac{1+d}{2}]\subseteq [1,d]$, which is non-zero in $H_0$. Thus,
$$\mathcal{N}\big(\ell(D_0)\big)=\mathcal{N}(D_0-I_0)=\{0\},\quad \mathcal{N}\big(\ell(D_1)\big)=\mathcal{N}(D_1-I_0)\ne\{0\}.
$$

We now prove that $\widetilde{Q_r}$ and $ q_r$ are not unitarily equivalent. By \eqref{equ:W4norm computation},
it suffices to show that the operators $A$ and $B$ are not unitarily equivalent, where
\begin{align*}A=\left(
                        \begin{array}{cc}
                          I_0 & 2\ell(D_0) \\
                          0 & 0 \\
                        \end{array}
                      \right),\quad B=\left(
                        \begin{array}{cc}
                          I_0 & 2\ell(D_{1}) \\
                          0 & 0 \\
                        \end{array}
                      \right).
\end{align*}
Assume for contradiction that there exists a unitary $U=\big(U_{ij}\big)_{1\le i,j\le 2}\in\mathbb{B}(H_0\oplus H_0)$ such that $UA=BU$.
A direct computation yields
\begin{equation}\label{equ:conditionofUij}U_{21}=0,\quad 2U_{11}\ell(D_0)=U_{12}+2\ell(D_{1})U_{22} .\end{equation}
Substituting  $U_{21}=0$ into the equations $UU^*=U^*U=I_0\oplus I_0$ yields
\begin{align*}
&U_{22}U_{22}^*=U_{11}U_{11}^*+U_{12}U^*_{12}=I_0,\quad U_{12}U_{22}^*=0,\\
&U_{11}^*U_{11}=U^*_{22}U_{22}+U_{12}^*U_{12}=I_0,\quad U_{11}^*U_{12}=0.\end{align*}
Left-multiplying both sides of the second equation in \eqref{equ:conditionofUij} by $U_{11}^*$ gives
$$\ell(D_0)=U_{11}^*\ell(D_{1})U_{22}.$$
Since $\mathcal{N}\big(\ell(D_0)\big)=\{0\}$, it follows that $\mathcal{N}(U_{22})=\{0\}$, so $\mathcal{N}(U_{22}^*U_{22})=\{0\}$.
As $U_{22}^*U_{22}$ is a  projection, this forces $U_{22}^*U_{22}=I_0$. Therefore,
$U_{22}$  is   unitary, which  implies $U_{12}=0$, and thus $U_{11}$ is also  unitary.
Consequently,
$$\ell(D_{1})=U_{11}\ell(D_0)U^*_{22}.$$
Since  $U_{11}$, $\ell(D_0)$ and $U^*_{22}$ are all  injective, it follows that $\mathcal{N}\big(\ell(D_{1})\big)=\{0\}$,
which is a contradiction.
\end{proof}

\section{New characterizations of the numerical ranges}\label{sec:numerical range}

\subsection{An operator version of the elliptical range theorem}\label{sub:operator version}

The elliptical range theorem is originally obtained in the matrix case, which is a powerful tool in dealing with the numerical ranges of matrices. The main purpose of this subsection is to set up an operator version of the elliptical range theorem.

The support function plays a crucial role in convex analysis \cite{Rockafellar,Schneider} within Euclidean space $\mathbb{R}^n$. By identifying $\mathbb{C}$  with $\mathbb{R}^2$, the support function for a convex subset $E$ of $\mathbb{C}$ is defined as follows.

\begin{definition}\label{defn of sf}{\rm \cite[Definition~2.4]{Klaja}}
 Let  $E$ be a bounded convex subset of $\mathbb{C}$. The support function of $E$, denoted by $h_E$, is defined on $\mathbb{R}$ by
\begin{equation*}
  h_E(\alpha)=\sup\{\mbox{Re}(ze^{-i\alpha}):z\in E\},\quad \alpha\in\mathbb{R}.
\end{equation*}
\end{definition}

From the above definition, it follows immediately that
\begin{equation}\label{equ:supportequal01}h_E(\alpha)=h_{\overline{E}}(\alpha),\quad \alpha\in\mathbb{R}.
\end{equation}
Moreover, if $\{E_i\}_{i\in I}$ is a family of uniformly bounded convex subsets of $\mathbb{C}$, then  it is readily verified that
\begin{equation}\label{equ:supportequal02}h_{\mbox{conv}\big(\cup_{i\in I} E_i\big)}(\alpha)=\sup_{i\in I} h_{E_i}(\alpha), \quad \alpha\in \mathbb{R}.
\end{equation}

The following lemma shows that the support function characterizes the closure of a convex set, which can be derived by a straightforward application of \cite[Corollary~1.2.12]{KR}.
\begin{lemma}\label{lem:in closure}{\rm \cite[Proposition~2.5]{Klaja}}
Let $E$ be a bounded convex subset of $\mathbb{C}$. Then
$$\overline{E}=\big\{z\in \mathbb{C}: \mbox{Re}(ze^{-i\alpha})\le h_E(\alpha),\quad  \forall \alpha\in \mathbb{R}\big\}.$$
\end{lemma}

\begin{remark}\label{rem:supportequal03}{\rm   Let $E_1$ and $E_2$ be two bounded convex subsets of $\mathbb{C}$.
A direct application of \eqref{equ:supportequal01} and Lemma~\ref{lem:in closure} yields
  \begin{equation}\label{equ:supportequal05}
  \overline{E_1}=\overline{E_2}\Longleftrightarrow h_{E_1}(\alpha)=h_{E_2}(\alpha),\quad  \forall \alpha\in \mathbb{R}.
\end{equation}
By virtue of this characterization and \eqref{equ:supportequal02}, we conclude that
for any bounded convex subset $S$ of $\mathbb{C}$ and any collection $\{E_i\}_{i\in I}$ of  uniformly bounded convex subsets of $\mathbb{C}$,
\begin{equation}\label{equ:supportequal04}
  \overline{S}=\overline{\mbox{conv}\big(\cup_{i\in I}E_i\big)}\Longleftrightarrow h_S(\alpha)=\sup_{i\in I}h_{E_i}(\alpha),\quad  \forall\alpha\in \mathbb{R}.
\end{equation}
}\end{remark}
Let $E$ be an elliptical disk whose major axis is parallel to one of the coordinate axes. In this case,
\begin{equation}\label{equ:elliptical disk001}
  x+iy\in E\Longleftrightarrow \frac{(x-x_0)^2}{a^2}+\frac{(y-y_0)^2}{b^2}\leq1,
\end{equation}
where  $x_0,y_0\in \mathbb{R}$ and $a,b\in (0,+\infty)$ are fixed parameters, while $x,y\in\mathbb{R}$ are variables. It is well-known that the support function of $E$ is given by
\begin{equation}\label{equ:sf of ellipse}
  h_E(\alpha)=
  x_0\cos\alpha+y_0\sin\alpha+\sqrt{a^2\cos^2\alpha+b^2\sin^2\alpha},\quad  \alpha\in [0,2\pi).
\end{equation}

For completeness, we state the following lemma, which identifies for each $\alpha\in [0,2\pi)$ a point
 $z_{\alpha}\in E$
    such that $
  h_E(\alpha)=\text{Re}(z_{\alpha}e^{-i\alpha})$.

\begin{lemma}\label{lem:sf of elliptical disk}
Let $E$ be the elliptical disk given  by \eqref{equ:elliptical disk001}. Then
 for every $\alpha\in [0,2\pi)$, equality
\eqref{equ:sf of ellipse} holds with $
  h_E(\alpha)=\text{Re}(z_{\alpha}e^{-i\alpha})$, where
\begin{equation}\label{equ:zalpha}
  z_{\alpha}=x_0+iy_0+\frac{a^2\cos\alpha+ib^2\sin\alpha }{\sqrt{a^2\cos^2\alpha+b^2\sin^2\alpha}}.
\end{equation}
\end{lemma}
\begin{proof}Elements of $E$ can be parameterized as $z=x(\theta,r)+iy(\theta,r)$, where
\begin{equation*}
\begin{cases}
x(\theta,r)=x_0+ar\cos \theta,  \\
y(\theta,r)=y_0+br\sin \theta,
\end{cases}\quad  r\in [0,1], \theta \in [0,2\pi).
\end{equation*}
 For each $\alpha\in [0,2\pi)$, by Definition~\ref{defn of sf} we have
\begin{align*}
  h_E(\alpha)=\sup\limits_{\theta,r }\left\{{x(\theta,r)\cos \alpha+y(\theta,r)\sin\alpha}\right\}=\sup\limits_{\theta,r }g(\theta,r),
\end{align*}
 where
 $$g(\theta,r)=x_0\cos\alpha+y_0\sin\alpha+r(a\cos\alpha\cos\theta+b\sin\alpha\sin\theta).$$
 Choose $\theta_0\in [0,2\pi)$ such that
$$\cos \theta_0=\frac{a\cos\alpha}{\sqrt{a^2\cos^2\alpha+b^2\sin^2\alpha}},\quad \sin \theta_0=\frac{b\sin\alpha}{\sqrt{a^2\cos^2\alpha+b^2\sin^2\alpha}}.$$
Then it follows that
\begin{align*}
  r(a\cos\alpha\cos\theta+b\sin\alpha\sin\theta) & =r\sqrt{a^2\cos^2\alpha+b^2\sin^2\alpha}\cdot\cos(\theta-\theta_0), \\
   & \le \sqrt{a^2\cos^2\alpha+b^2\sin^2\alpha}.
\end{align*}
Moreover, equality holds if and only if $r=1$ and $\theta=\theta_0$.
Therefore,
\begin{align*}
\sup\limits_{\theta,r}g(\theta,r)=\text{Re}(z_{\alpha}e^{-i\alpha})=x_0\cos\alpha+y_0\sin\alpha+\sqrt{a^2\cos^2\alpha+b^2\sin^2\alpha}.
\end{align*}
where $z_{\alpha}$ is defined in \eqref{equ:zalpha}.
This completes the proof.
\end{proof}

\begin{remark} Several known results (e.g.,\,\cite[Lemma~3.4]{Klaja}), follow directly from \eqref{equ:sf of ellipse}.
\end{remark}

As a  complement to Lemma~\ref{lem:sf of elliptical disk}, the following lemma addresses the support functions of line segments. These segments can be regarded as
degenerate cases of elliptical disks,  occurring when  either $b=0$ or $a=0$.

\begin{lemma}\label{lem:sf of interval}
Let $E$ and $F$ be the line segments  defined  by
\begin{align*}&E=\{x+iy_0:x\in\mathbb{R}, x_0-a\le x\le x_0+a\},\\
&F=\{x_0+iy:y\in\mathbb{R}, y_0-b\le y\le y_0+b\},
\end{align*}
where $x_0,y_0\in \mathbb{R}$ and $a,b\in (0,+\infty)$ are fixed. Then the support functions $h_E(\alpha)$ and  $h_F(\alpha)$  correspond to the degenerate forms  of  \eqref{equ:sf of ellipse} with $b=0$ and $a=0$, respectively.
\end{lemma}
\begin{proof}
The result follows directly from the definition of the support function and the parametric descriptions of $E$ and $F$.
\end{proof}

Recall that for every $T\in \mathbb{B}(H)$, its numerical range $W(T)$ and numerical radius $w(T)$ are defined respectively by
\begin{equation*}
  W(T) = \big\{\left<Tx,x\right>:x\in H, \|x\| = 1\big\},\quad w(T)=\sup\big\{|\lambda|:\lambda\in W(T)\big\}.
\end{equation*}
Note that $w(T)=\sup\{|\lambda|:\lambda\in \overline{W(T)}\}$.
By the classical Toeplitz-Hausdorff theorem \cite[Theorem~1.1-2]{GR}, the numerical range $W(T)$ is convex for every $T\in\mathbb{B}(H)$. Consequently, for any $T_i\in\mathbb{B}(H_i)$ $(i=1,2)$,
\begin{equation}\label{equ:equal W}
W(T_1\oplus T_2)=\mbox{conv}\big(W(T_1) \cup W(T_2)\big).
\end{equation}

For the case of matrices, we now state the so-called elliptical range theorem. A proof can be found  in \cite[Theorem]{Li}.

\begin{lemma}\label{lem:the numeriacal of matrix01}
Let $A$ be a $2\times2$ matrix with eigenvalues $\lambda_1$ and $\lambda_2$. Then the  numerical range of $A$ is an elliptical disk whose foci are $\lambda_1$ and $\lambda_2$, and whose minor axis has length $$\big(\text{tr}(A^*A)-|\lambda_1|^2-|\lambda_2|^2\big)^{\frac12},$$ where $\text{tr}(A^*A)$ denotes the trace of $A^*A$.
\end{lemma}

In some cases, the numerical range of a matrix is contained in that of an operator. We provide  such a result below.

\begin{lemma}\label{lemma:point-spectrum-numerical-range}
Let $D \in \mathbb{B}(K)_{\text{sa}}$ and let
$f=(f_{ij})_{1\leq i,j\leq 2}\in M_2\Big(C\big(\sigma(D)\big)\Big)$. If $t_0 \in \sigma_p(D)$, then $W(f_{t_0}) \subseteq W\big(f(D)\big)$,
where  $f_{t_0}$ and $f(D)$ are defined by \eqref{2ff}.
\end{lemma}
\begin{proof}
Suppose  $t_0 \in \sigma_p(D)$. Then there exists a unit vector $u \in K$ such that $Du = t_0 u$,
which gives  $p(D)u=p(t_0)u$ for every polynomial $p$.
 By continuous functional calculus, it follows that
\begin{equation}\label{equ:continuous functional calculus}
g(D)u=g(t_0)u, \quad  \forall  g\in C\big(\sigma(D)\big).
\end{equation}

Now let $z\in W(f_{t_0})$. There exist $\xi,\eta\in \mathbb{C} $ with $|\xi|^2+|\eta|^2=1$ such that
$$z
=\left\langle f_{t_0}(\xi,\eta)^T,(\xi,\eta)^T\right\rangle.$$
Define $\xi'=\xi u$ and $\eta'=\eta u$. Then $\xi',\eta'\in K$ and satisfy  $\|\xi'\|^2+\|\eta'\|^2=1$. Using  \eqref{2ff} and  \eqref{equ:continuous functional calculus}, we compute:
\begin{align*}
\left\langle f(D)(\xi',\eta' )^T,(\xi',\eta' )^T\right\rangle=&
\left\langle \binom{\xi f_{11}(D)u+\eta f_{12}(D)u}{\xi f_{21}(D)u+\eta f_{22}(D)u},\binom{\xi u}{\eta u}\right\rangle\\
=&\left\langle \binom{\xi f_{11}(t_0)u+\eta f_{12}(t_0)u}{\xi f_{21}(t_0)u+\eta f_{22}(t_0)u},\binom{\xi u}{\eta u}\right\rangle\\
=&
f_{11}(t_0)|\xi|^2+f_{12}(t_0)\bar{\xi}\eta+f_{21}(t_0)\xi\bar{\eta}+f_{22}(t_0)|\eta|^2\\
=&
\left\langle f_{t_0}(\xi,\eta)^T,(\xi,\eta)^T\right\rangle=z.
\end{align*}
So, $z\in W\big(f(D)\big)$.
By the arbitrariness of $z$, we get  $W(f_{t_0})\subseteq   W\big(f(D)\big)$.
\end{proof}

We now present the following lemma, which is essential for establishing the main result of this subsection.

\begin{lemma}\label{lem:the norm of T+}
Let $T \in \mathbb{B}(H)_{\mathrm{sa}}$ and suppose that its positive part $T_{+}$ is nonzero. Then
\[
\sup\bigl\{\langle Tx, x\rangle : x \in H,\ \|x\| = 1\bigr\} = \|T_{+}\|.
\]
\end{lemma}

\begin{proof}
Decompose $T = T_{+} - T_{-}$ with $T_{+}T_{-} = T_{-}T_{+} = 0$, where $T_{+}$ and $T_{-}$ are the positive and negative parts of $T$, respectively. Define
\[
\beta = \sup\bigl\{\langle Tx, x\rangle : x \in H,\ \|x\| = 1\bigr\}.
\]
For any unit vector $x \in H$, we have $\langle Tx, x\rangle \le \langle T_{+}x, x\rangle \le \langle \|T_{+}\|x, x\rangle = \|T_{+}\|$. Hence $\beta \le \|T_{+}\|$.

To establish the reverse inequality $\|T_{+}\| \le \beta$, note that since $T_{+} \neq 0$,
\begin{align}\label{equ:the norm of T+}
\|T_{+}\| &= w(T_{+}) = \sup\bigl\{\langle T_{+}x, x\rangle : x \in H,\ \|x\| = 1\bigr\}\nonumber \\
&= \sup\bigl\{\langle T_{+}x, x\rangle : x \in \overline{\mathcal{R}(T_{+})},\ \|x\| = 1\bigr\}.
\end{align}
Now, for any unit vector $x \in \overline{\mathcal{R}(T_{+})}$, we have $T_{-}x = 0$, and therefore
\[
\langle T_{+}x, x\rangle = \langle Tx, x\rangle \le \beta.
\]
Combining this with \eqref{equ:the norm of T+} yields $\|T_{+}\| \le \beta$. Consequently, $\beta = \|T_{+}\|$.
\end{proof}

\begin{remark}
If $T \in \mathbb{B}(H)$ is positive, then  Lemma~\ref{lem:the norm of T+} yields
$w(T)=\|T\|$, where $w(T)$ denotes the numerical radius of $T$.
\end{remark}

We are now ready to present an operator version of the elliptical range theorem.

\begin{theorem}\label{thm:the numerical of matrix02}
Let $D \in \mathbb{B}(K)_{\mathrm{sa}}$ and let $f = (f_{ij})_{1\le i,j\le 2} \in M_2\Big(C\big(\sigma(D)\big)\Big)$. Then
\begin{equation}\label{equ:the num of B}
  \overline{W\big(f(D)\big)} = \overline{\operatorname{conv}\Big( \bigcup_{t\in \sigma(D)} W(f_t) \Big)},
\end{equation}
where $f(D)$ and $f_t$ for $t \in \sigma(D)$ are defined as in \eqref{2ff}.
\end{theorem}

\begin{proof}
We begin by recalling some basic facts about $C^*$-algebras. Let $\mathfrak{A}$ and $\mathfrak{B}$ be $C^*$-algebras and let $\pi : \mathfrak{A} \to \mathfrak{B}$ be a $C^*$-homomorphism. For any $a \in \mathfrak{A}$, we have $\pi(|a|) = |\pi(a)|$. In particular, if $a = a^*$, then
\[
\pi(a_{+}) = \pi\!\left(\frac{a + |a|}{2}\right) = \frac{\pi(a) + |\pi(a)|}{2} = \big(\pi(a)\big)_{+}.
\]

Now, let $\mathfrak{A} = M_2\big(C^*\{D, I_K\}\big)$, which is a $C^*$-subalgebra of $\mathbb{B}(K \oplus K)$. Since $D$ is self-adjoint, we can write
\[
\mathfrak{A} = \left\{
\begin{pmatrix}
g_{11}(D) & g_{12}(D) \\
g_{21}(D) & g_{22}(D)
\end{pmatrix}
: g_{ij} \in C\big(\sigma(D)\big),\; 1 \le i,j \le 2
\right\}.
\]
For $g(D) = \big(g_{ij}(D)\big) \in \mathfrak{A}$, it follows from the proof of Lemma~\ref{lem:technical lem} that
\[
\|g(D)\| = \max_{t \in \sigma(D)} \|g_t\|.
\]
Let $\mathfrak{B} = M_2(\mathbb{C})$. For each $t \in \sigma(D)$, define the $C^*$-homomorphism $\pi_t : \mathfrak{A} \to \mathfrak{B}$ by $\pi_t\big(g(D)\big) = g_t$. Then
\begin{equation}\label{equ:the norm of g}
\|g(D)\| = \max_{t \in \sigma(D)} \big\|\pi_t\big(g(D)\big)\big\|.
\end{equation}

Now let $f(D)$ and $f_t$ be as above. For any $\alpha \in [0, 2\pi)$ and any unit vector $\xi \in K \oplus K$, we have
\[
\left|\left\langle \operatorname{Re}\big(e^{-i\alpha}f(D)\big)\xi, \xi \right\rangle\right|
\le \big\|\operatorname{Re}\big(e^{-i\alpha}f(D)\big)\big\|
\le \|e^{-i\alpha}f(D)\| = \|f(D)\|.
\]
Define $g(D) = f(D) + \lambda_\alpha I_{K\oplus K}$, where
\[
\lambda_\alpha=
\begin{cases}
\dfrac{\|f(D)\|+1}{\cos\alpha}, & \text{if } \cos\alpha \ne 0, \\[8pt]
(\|f(D)\|+1)i, & \text{if } \alpha = \dfrac{\pi}{2}, \\[8pt]
-(\|f(D)\|+1)i, & \text{if } \alpha = \dfrac{3\pi}{2}.
\end{cases}
\]
Then
\begin{equation}\label{from f 2 g}
\operatorname{Re}\big(e^{-i\alpha}g(D)\big) = \operatorname{Re}\big(e^{-i\alpha}f(D)\big) + (\|f(D)\|+1)I_{K\oplus K}.
\end{equation}
Thus $\operatorname{Re}\big(e^{-i\alpha}g(D)\big)$ is positive definite. Since each $\pi_t$ is a unital $C^*$-homomorphism, for every $t \in \sigma(D)$ we have
\begin{align*}
\operatorname{Re}\big(e^{-i\alpha}g_t\big)
&= \operatorname{Re}\!\Big[ e^{-i\alpha} \pi_t\big(g(D)\big) \Big]
 = \operatorname{Re}\!\Big[ \pi_t\big(e^{-i\alpha}g(D)\big) \Big] \\
&= \pi_t\!\Big[ \operatorname{Re}\big(e^{-i\alpha}g(D)\big) \Big].
\end{align*}
Consequently, $\operatorname{Re}\big(e^{-i\alpha}g_t\big)$ is also positive definite for all $t \in \sigma(D)$.

Fix $\alpha \in [0, 2\pi)$. By Definition~\ref{defn of sf}, Lemma~\ref{lem:the norm of T+}, and \eqref{equ:the norm of g},
\begin{align*}
h_{W\big(g(D)\big)}(\alpha)
&= \sup\left\{ \operatorname{Re}\big(e^{-i\alpha}\langle g(D)\xi, \xi\rangle\big) : \xi \in K \oplus K,\ \|\xi\| = 1 \right\} \\
&= \sup\left\{ \left\langle \operatorname{Re}\big(e^{-i\alpha}g(D)\big)\xi, \xi \right\rangle : \xi \in K \oplus K,\ \|\xi\| = 1 \right\} \\
&= \Big\| \operatorname{Re}\big(e^{-i\alpha}g(D)\big) \Big\|
 = \max_{t \in \sigma(D)} \Big\| \pi_t\!\left[ \operatorname{Re}\big(e^{-i\alpha}g(D)\big) \right] \Big\| \\
&= \max_{t \in \sigma(D)} \big\| \operatorname{Re}(e^{-i\alpha}g_t) \big\|.
\end{align*}
Similarly, for each $t \in \sigma(D)$,
\[
h_{W(g_t)}(\alpha) = \sup\left\{ \operatorname{Re}\big(e^{-i\alpha}\langle g_t\xi, \xi\rangle\big) : \xi \in \mathbb{C}^2,\ \|\xi\| = 1 \right\}
= \big\| \operatorname{Re}(e^{-i\alpha}g_t) \big\|.
\]
Hence,
\[
h_{W\big(g(D)\big)}(\alpha) = \max_{t \in \sigma(D)} h_{W(g_t)}(\alpha).
\]
It follows from \eqref{from f 2 g} that
\begin{align*}
&h_{W\big(f(D)\big)}(\alpha)= h_{W\big(g(D)\big)}(\alpha) - (\|f(D)\|+1), \\
&h_{W(f_t)}(\alpha)= h_{W(g_t)}(\alpha) - (\|f(D)\|+1).
\end{align*}
Therefore, for every $\alpha \in [0, 2\pi)$,
\[
h_{W\big(f(D)\big)}(\alpha) = \sup_{t \in \sigma(D)} h_{W(f_t)}(\alpha).
\]
The result now follows from \eqref{equ:supportequal04}.
\end{proof}

\begin{remark}
  By  Lemma~\ref{lem:the numeriacal of matrix01}, the numerical range  $W(f_t)$ is an elliptical disk  for each $t\in \sigma(D)$. Therefore, Theorem~\ref{thm:the numerical of matrix02}  implies that the closure of the numerical range of $f(D)$ coincides with the closure of the convex hull of the union of a family of elliptical disks
\end{remark}

\begin{remark}\label{rem:useful remark001}
Note that, in general, equality \eqref{equ:the num of B} does not hold without the closures, i.e.,
 \begin{equation*}
  W\big(f(D)\big)=\mbox{conv}\big(\cup_{t\in \sigma(D)}W(f_{t})\big)
\end{equation*}
is not always true. For example, as will be shown in Corollary~\ref{cor:num of Qr}, for each $r>1$, the numerical range $W(T_{{}_{Q_r}})$ is an open elliptical disk described by \eqref{equ:open ellipse-r}, whereas $\|D_0\|\in\sigma(D_0)$ and $W(f_{\|D_0\|})$ equals the closure of this elliptical disk.
\end{remark}

\subsection{The numerical ranges of certain non-quadratic operators }\label{subsec:non-quadratic}

To derive the main results of this subsection, we employ the lemma.
\begin{lemma}\label{lem:convtos}{\rm (cf.\,\cite[Theorem~1.1.11]{Schneider})}
 Let   $C$  be a  bounded subset of $\mathbb{C}$. Then
 $\overline{\text{conv}(C)}=\text{conv}(\overline{C})$.
\end{lemma}

In the following theorem, we introduce a non-quadratic operator $T_{{}_Q}$ by using  the matched projection of an  idempotent, and then characterize the closure of its numerical range.

\begin{theorem}\label{thm:the num of Q2}
Let $Q\in \mathbb{B}(H)$ be a non-projection idempotent. Define $T_{{}_Q}\in\mathbb{B}(H)$ by
\begin{equation}\label{equ:nonquadopt}T_{{}_Q}=m(Q)+m(Q)Q.\end{equation}
Then $T_{{}_Q}$ is not a quadratic operator. Moreover, $\|T_{{}_Q}\|=2a$, $w(T_{{}_Q})=a+c$, and $\overline{W(T_{{}_Q})}$ is a non-degenerate elliptical disk given by
\begin{equation}\label{equ:useful003}\overline{W(T_{{}_Q})}=\Big\{x+iy:x,y\in\mathbb{R}, \frac{(x-c)^2}{a^2}+\frac{y^2}{b^2}\leq1\Big\},\end{equation}
where
\begin{equation}\label{equ:abcinE}d=\frac{\|Q\|+1}{2},\quad
  a=\frac{\sqrt{2d^2+d+1}}{2},\quad b=\frac{\ell(d)}{2},\quad c=\frac{d+1}{2}.
\end{equation}
\end{theorem}
\begin{proof}Let $\widetilde{Q}$ be defined as in  \eqref{eq:new decomposition of Q2 hat2}, where the operator  $D\in \mathbb{B}(H_5)$ satisfies  \eqref{equ:2propofD}.
Set $d=\|D\|$ and define $T_{{}_Q}$ as above. Then $d=\frac12(\|Q\|+1)>1$. From  \eqref{equ:the unitarily of M(Q)}--\eqref{eq:new decomposition of Q2 hat2}, we have
$$T_{{}_Q}=V_Q^*\left[2I_{H_1}\oplus 0\oplus \widetilde{T_{{}_Q}}\right]V_Q,$$
where $V_Q$ is a unitary operator and
\begin{equation}\label{equ:tildeT}
\widetilde{T_{{}_Q}}:=f(D)=\left(
           \begin{array}{cc}
             D+I_{H_5} & -\ell(D) \\
             0 & 0 \\
           \end{array}
         \right).
\end{equation}
Since $D\ge I_{H_5}$, for any $a,b\in\mathbb{C}$, a direct computation shows that
$$(T_{{}_Q}-aI_H)(T_{{}_Q}-bI_H)=0\Longleftrightarrow ab=0,\quad  a+b=2,\quad  D=I_{H_5}.$$
Note that the condition  $D=I_{H_5}$ holds only if  $Q$ is a projection, so $T_{{}_Q}$ is not a quadratic operator.

From the expressions obtained for $T_{{}_Q}$ and $\widetilde{T_{{}_Q}}$, we have
\begin{align*}\|T_{{}_Q}\|&=\max\Big\{\|2I_{H_1}\|,\sqrt{\big\|\widetilde{T_{{}_Q}}\widetilde{T_{{}_Q}}^*\big\|}\Big\}=\sqrt{2d^2+d+1}.\end{align*}
Moreover, a direct computation shows that for $\alpha\in [0,2\pi)$,
\begin{align*}\mbox{Re}\big(e^{-i\alpha}\widetilde{T_{{}_Q}}\big)=\begin{pmatrix}\cos\alpha\cdot (D+I_{H_5})&-\frac12 e^{-\alpha i}\ell(D)\\
-\frac12 e^{\alpha i}\ell(D)&0\end{pmatrix}.
\end{align*}
Suppose $\mbox{Re}\big(e^{-i\alpha}\widetilde{T_{{}_Q}}\big)\le 0$. Then for all $x,y\in H_5$,
$$\big\langle \mbox{Re}\big(e^{-i\alpha}\widetilde{T_{{}_Q}}\big) (x,y)^T,(x,y)^T\big\rangle\le 0,$$
or equivalently,
$$\cos\alpha\cdot \langle (D+I_{H_5})x,x\rangle-\frac12\mbox{Re}\big(e^{-\alpha i}\langle \ell(D)y,x\rangle\big)\le 0.$$
Replacing $y$ by $ky$ for any $k\in\mathbb{R}$ forces
$$\mbox{Re}\big(e^{-\alpha i}\langle \ell(D)y,x\rangle\big)=0\quad \mbox{and}\quad \cos\alpha\cdot \langle (D+I_{H_5})x,x\rangle=0$$
 for all $x,y\in H_5$.
Since $D+I_{H_5}$ is positive definite,  we conclude that
$\cos\alpha=0$. It then follows that $\mbox{Re}\langle \ell(D)(\pm i y),x\rangle=0$ for all $x,y\in H_5$, which implies
$\ell(D)=0$. This contradicts the fact that $\|\ell(D)\|=\ell(d)>0$. Therefore, for every $\alpha\in [0,2\pi)$,
the positive part of $\mbox{Re}\big(e^{-i\alpha}\widetilde{T_{{}_Q}}\big)$ is nonzero.

By \eqref{equ:equal W}, \eqref{equ:tildeT} and  Theorem~\ref{thm:the numerical of matrix02}, we have
\begin{align}
\label{equ:the num of mq001}
&W(T_{{}_Q})=\mbox{conv}\big[W(2I_{H_1}\oplus 0) \cup W(\widetilde{T_{{}_Q}})\big],\\
\label{equ:the num of mq002}&\overline{W(\widetilde{T_{{}_Q}})}=\overline{\mbox{conv}\big(\cup_{t\in \sigma(D)}W(f_{t})\big)},
\end{align}
where $\ell(t)$ is defined by \eqref{defn of lt}, and $f_t\in M_2(\mathbb{C})$ is given by
$$f_t=\left(
        \begin{array}{cc}
          t+1 & -\ell(t) \\
          0 & 0 \\
        \end{array}
      \right),\quad  t\in \sigma(D).$$
A simple computation yields
 $$\text{tr}(f_{t}f_{t}^*)-|\lambda_1|^2-|\lambda_2|^2=\ell^2(t),$$
in which $\lambda_1=0$ and $\lambda_2=t+1$ are the eigenvalues of $f_t$.
Hence, by  Lemma~\ref{lem:the numeriacal of matrix01} and \eqref{equ:sf of ellipse},  $W(f_t)$ is an elliptical disk (including the
degenerate case when $t=1\in\sigma(D)$) such that for any $t\in\sigma(D)$ and $\alpha\in [0,2\pi)$,
$$h_{W(f_t)}(\alpha)=x_t\cos\alpha+\sqrt{a_t^2 \cos^2\alpha+b_t^2\sin^2\alpha},
$$
with
\begin{align*}
  x_t= & c_t=\frac{1}{2}(t+1),\quad   b_t = \frac{1}{2}\ell(t),\quad   a_t=\sqrt{b_t^2+c_t^2}=\frac12\sqrt{2t^2+t+1}.
\end{align*}
It follows that
\begin{equation*}
  h_{W(f_t)}(\alpha)=\frac{t+1}{2}\cos\alpha+\frac{1}{2}\sqrt{(t+1)^2\cos^2\alpha+\ell^2(t)}.
\end{equation*}

Suppose $\cos\alpha\ne 0$. By direct  computation, we obtain
\begin{align*}
  \frac{\partial [h_{W(f_t)}(\alpha)]}{\partial t} &= \frac{g(t,\alpha)}{2\sqrt{(t+1)^2\cos^2\alpha+\ell^2(t)}},
\end{align*}
where \begin{align*}g(t,\alpha)&=(t+1)\cos^2\alpha+t-\frac12+\sqrt{(t+1)^2\cos^2\alpha+\ell^2(t)}\cos\alpha\\
&\ge (t+1)\cos^2\alpha+t-\frac12-\sqrt{(t+1)^2\cos^2\alpha+\ell^2(t)}\sqrt{\cos^2\alpha}\\
&\ge (t+1)\cos^2\alpha+\ell(t)-\sqrt{(t+1)^2\cos^4\alpha+\ell^2(t)}.
\end{align*}
Since $t\ge 1$, it follows that $g(t,\alpha)\ge 0$.
Therefore,
$$ h_{W(f_{t_1})}(\alpha) \le h_{W(f_{t_2})}(\alpha), \quad \forall t_1,t_2\in \sigma(D)\ \mbox{with}\  t_1\le t_2.$$
If $\cos\alpha=0$, then $h_{W(f_t)}(\alpha)=\ell(t)/2$, so the above inequality holds trivially.
By Lemma~\ref{lem:in closure}, this implies $W(f_{t_1})\subseteq \overline{W(f_{t_2})}=W(f_{t_2})$ for such $t_1$ and $t_2$.
In view of $\max\limits_{t\in \sigma(D)}t=d$,  we have
\begin{equation*}\bigcup_{t\in \sigma(D)}W(f_{t})=W(f_{d}),
\end{equation*}
which is the elliptical disk given by the right-hand side of \eqref{equ:useful003}, with parameters $a,b$ and $c$ defined in \eqref{equ:abcinE}.

Combining \eqref{equ:the num of mq002} and \eqref{equ:useful003} yields
$$\overline{W\big(\widetilde{T_{{}_Q}}\big)}=\overline{\mbox{conv}\big(W(f_d)\big)}=W(f_d).$$
As $d>1$, it follows from \eqref{equ:useful003} that
$0,2\in \text{int}\big(W(f_d)\big)$.
Note that  \begin{equation*}\label{equ:W(I)=0}
 W\big(I_{H_1}\big)= \begin{cases}
  1, & \mbox{if $H_1\ne 0$, }\\
  0, & \mbox{otherwise}.
\end{cases}
\end{equation*}
Therefore, $W(2I_{H_1}\oplus 0)\subseteq [0,2]\subseteq \text{int}\big(W(f_d)\big)\subseteq W(f_d)$.
It can be deduced from \eqref{equ:the num of mq001} and Lemma~\ref{lem:convtos} that
\begin{align*}
\overline{W(T_{{}_Q})}=&\overline{\mbox{conv}\big[W\big(2I_{H_1}\oplus 0\big)\cup W(\widetilde{T_{{}_Q}})\big]}=\mbox{conv}\big[\overline{W\big(2I_{H_1}\oplus 0
\big)}\cup \overline{W(\widetilde{T_{{}_Q}})}\big]\\
=&\mbox{conv}\big[\overline{W\big(2I_{H_1}\oplus 0
\big)}\cup W(f_d)\big]=\text{conv}\big(W(f_d)\big)=W(f_d).
\end{align*}

Geometrically, the numerical radius of $T_{{}_Q}$ is given by
$w(T_{{}_Q})=c+a$. This completes the proof.
\end{proof}

To prove $\mbox{int}\big(\overline{W(T_{{}_Q})}\big)\subseteq W(T_{{}_Q})$,  we invoke the following lemma initially established in \cite{WX}.

\begin{lemma}\label{lem:interior contained F} Let  $E$ be a non-degenerate elliptical disk and $F$ be a convex subset of $E$. If $F$ is dense in $E$ and contains the center of $E$, then the interior of $E$ is contained in $F$.
\end{lemma}

\begin{tikzpicture}[scale=1.4, rotate=20, >=stealth]

    \draw[thick] (0,0) ellipse (3.5 and 2.2);

    \coordinate (O) at (0,0);
    \fill (O) circle (1pt) node[below left] {$O$};

    \coordinate (z) at (1.2,0.8);
    \coordinate (z') at (2.4,1.6);

    \draw[dashed] (O) -- (z');
    \fill (z) circle (0.7pt) node[below] {$z$};
    \fill (z') circle (0.7pt) node[above right] {$z'$};

    \draw[dashed, blue, thick] (z) circle (1.2);
    \node[blue, above right, font=\scriptsize] at (-0.4,1.8) {$B(z,r)$};

    \coordinate (A) at ($(z)!3cm!90:(O)$);
    \coordinate (B) at ($(z)!3cm!-90:(O)$);
    \draw[thick, gray] (A) -- (B);

    \node[font=\small] at ($(z)!0.5!135:(O)$) {(i)};
    \node[font=\small] at ($(z)!0.6!225:(O)$) {(ii)};
    \node[font=\small] at ($(z)!0.5!315:(O)$) {(iii)};
    \node[font=\small] at ($(z)!0.5!45:(O)$) {(iv)};

    \coordinate (u1) at (1.0, 1.65);
    \coordinate (u2) at (2.1, 1);
    \fill[purple] (u1) circle (0.7pt) node[below right] {$u_1$};
    \fill[purple] (u2) circle (0.7pt) node[below right] {$u_2$};
    \draw[purple, thick] (u1) -- (u2);

    \coordinate (u3) at ($(u1)!0.71!(u2)$);
    \fill[purple] (u3) circle (0.7pt) node[above] {$u_3$};

\end{tikzpicture}

\begin{proof}
As illustrated in the figure, let $O$ be the center of the elliptical disk $E$, and take an arbitrary point $z \in \operatorname{int}(E)$. The ray $Oz$ meets the boundary $\partial E$ at $z'$. Since $z$ is an interior point, there exists $r>0$ such that $B(z,r)\subseteq E$.

Through $z$, draw a line perpendicular to $Oz$. This line divides $B(z,r)$ into four regions: (i), (ii), (iii), and (iv), as illustrated. By the density of $F$ in $E$, we may choose $u_1 \in F$ from region (i) and $u_2 \in F$ from region (ii). The segment $u_1u_2$ intersects $zz'$ at a point $u_3$. By convexity, $u_3 \in F$. Since $O \in F$ and $F$ is convex, the whole segment $Ou_3$ lies in $F$. As $z$ lies between $O$ and $u_3$, it follows that $z \in F$.
Because $z$ was arbitrary, we conclude $\operatorname{int}(E) \subseteq F$.
\end{proof}

\begin{theorem}\label{thm:biggerintT}
Let  $Q\in \mathbb{B}(H)$ be a non-projection idempotent. Then
$$\mbox{int}\big(\overline{W(T_{{}_Q})}\big)\subseteq W(\widetilde{T_{{}_Q}})=W(T_{{}_Q}),$$ where
$T_{{}_Q}$ and $\widetilde{T_{{}_Q}}$  are defined by  \eqref{equ:nonquadopt}
and \eqref{equ:tildeT}, respectively.
\end{theorem}
\begin{proof}
Let   $T_{{}_Q}$, $\widetilde{T_{{}_Q}}$, $D$ and $d$ be defined as in the proof of Theorem~\ref{thm:the num of Q2}.
For any vectors $u,v\in H_5$,  Equation \eqref{equ:tildeT} gives
\begin{equation}\label{equ:innertilTTT}
\left\langle \widetilde{T_{{}_Q}}\binom{u}{v},\binom{u}{v}\right\rangle=\big\langle (D+I_{H_5})u,u\big\rangle-\big\langle \ell(D)v,u\big\rangle.
\end{equation}

We now prove that the center $\frac{d+1}{2}$ of $\overline{W(T_{{}_Q})}$ lies in $W(\widetilde{T_{{}_Q}})$.
 Let $y\in H_5$ be any unit vector. Define $\zeta\in H_5\oplus H_5$ by $\zeta=(0,y)^T$. Then by \eqref{equ:innertilTTT},
 $\langle \widetilde{T_{{}_Q}}\zeta,\zeta\rangle=0$,
which implies  $0\in W(\widetilde{T_{{}_Q}})$.
Since $D+I_{H_5}$  is a positive operator with norm $d+1$, there exists  a unit vector $x_0\in H_5$ such that
$$t_0:=\big\langle (D+I_{H_5})x_0,x_0\big\rangle\in (d+3/4,d+1].$$
Let $\xi=(x_0,0)^T$. Then by \eqref{equ:innertilTTT}, $t_0=\langle \widetilde{T_{{}_Q}}\xi,\xi\rangle$, so
$t_0\in W(\widetilde{T_{{}_Q}})$.
Since $\frac{d+1}{2}\in [0,t_0]$ and $W(\widetilde{T_{{}_Q}})$ is convex, $\frac{d+1}{2}\in W(\widetilde{T_{{}_Q}})$.

From the proof of Theorem~\ref{thm:the num of Q2}, we have
$$\overline{W\big(\widetilde{T_{{}_Q}}\big)}=W(f_d)=\overline{W(T_{{}_Q})},$$
meaning  $W\big(\widetilde{T_{{}_Q}}\big)$ is dense in $\overline{W(T_{{}_Q})}$. It follows from Lemma~\ref{lem:interior contained F}  that $\mbox{int}\big(\overline{W(T_{{}_Q})}\big)\subseteq W(\widetilde{T_{{}_Q}})$.

The same proof of  Theorem~\ref{thm:the num of Q2} also shows $W(2I_{H_1}\oplus 0)\subseteq \mbox{int}\big(\overline{W(T_{{}_Q})}\big)$. Thus,
  $W(2I_{H_1}\oplus 0)\subseteq W(\widetilde{T_{{}_Q}})$.
  The equality $W(T_{{}_Q})=W(\widetilde{T_{{}_Q}})$ follows immediately from \eqref{equ:the num of mq001}.
\end{proof}

When $H$ is a  finite-dimensional Hilbert space, $W(A)$ is closed in $\mathbb{C}$ for all $A\in \mathbb{B}(H)$ \cite[Proposition~1.1]{WG}. However, when $H$ is infinite-dimensional,  $W(A)$ may fail to be closed for some  $A\in \mathbb{B}(H)$. This motivates us to investigate the closedness of $W(T_{{}_Q})$  using an approach distinct from that employed in \cite{TW,WX}. For this purpose, we require the following two lemmas.

\begin{lemma}\label{lem:Du=du001}
 Let   $D \in \mathbb{B}(K)$ satisfy $ D\ge I_K$ and $\|D\|>1$. Define $\ell(D)$ as in \eqref{equ:newtemp-003}.   Then for any $ u\in K$, $Du=\|D\|u$ if and only if $\ell(D)u=\ell(\|D\|)u$.
\end{lemma}
\begin{proof}Denote by $d=\|D\|$, and assume $u\in K$. Suppose first that $\ell(D)u=\ell(d)u$.
Then for any polynomial $p$, we have
$p\big(\ell(D)\big)u=p\big(\ell(d)\big)u$. By the continuous functional calculus, this implies
\begin{equation*}
f\big(\ell(D)\big)u=f\big(\ell(d)\big)u, \quad \forall f\in C\Big(\sigma\big(\ell(D)\big)\Big).
\end{equation*}
Define the function $$f(t)=\frac12\big(1+\sqrt{1+4t^2}\big),\quad  t\in\sigma\big(\ell(D)\big).$$ Then
 \begin{align*}Du=&\frac12\left[I_{K}+\big(I_{K}+4D(D-I_{K})\big)^\frac12\right]u=f\big(\ell(D)\big)u=f\big(\ell(d)\big)u\\
 =&
  \frac12\left[1+\big(1+4d(d-1)\big)^\frac12\right]u=du.
\end{align*}

Conversely, suppose $Du=du$. Then
\begin{equation*}\label{equ:eigenvalueequ}
g(D)u=g(d)u, \quad \forall g\in C\big(\sigma(D)\big).
\end{equation*}
 Since the function $t\to \ell(t)$ is continuous on $\sigma(D)$, we have
$\ell(D)u=\ell(d)u$. This completes the proof.
\end{proof}

\begin{lemma}\label{lem:Da and SD} Let   $D \in \mathbb{B}(K)$ satisfy $ D\ge I_K$ and $d:=\|D\|>1$. Define $\ell(D)$ and $\ell(d)$ as in \eqref{equ:newtemp-003},
and define  $D_a\in \mathbb{B}(K)$ by
\begin{align}\label{equ:Da and SD}
 D_a =\frac{D-aI_K}{\|D-aI_K\|},\quad  -\infty< a<\frac{d}{2d-1}.
\end{align}
Then $D_a$ is positive definite, and
\begin{equation*}
  D_a-\frac{\ell(D)}{\ell(d)}=A(dI_K-D),
\end{equation*}
where $A \in\mathbb{B}(K)$ is a positive definite operator given by
\begin{align*}\label{equ:A001}
  A = \Big[(d-a)^2\ell^2(d)\Big(D_a+\frac{\ell(D)}{\ell(d)}\Big)\Big]^{-1} \Big[(a^2-2da+d)D+a^2(d-1)I_K\Big].
\end{align*}
\end{lemma}
\begin{proof}
Suppose that $a\in \big(-\infty,\frac{d}{2d-1}\big)$. Define the function
\begin{equation*}\label{equ:funtionf01}
f(t)=\frac{t-a}{d-a}-\frac{\ell(t)}{\ell(d)},\quad t\in [1,d]. \end{equation*}
Then
\begin{align*}
  f(t) =\frac{(t-a)^2\ell^2(d)-(d-a)^2\ell^2(t)}{(d-a)\ell(d)\big[(t-a)\ell(d)+(d-a)\ell(t)\big]}=g(t)(d-t),
\end{align*}
where
$$g(t)=\frac{h(t)}
  {(d-a)^2\ell^2(d)\Big[\frac{t-a}{d-a}+\frac{\ell(t)}{\ell(d)}\Big]},$$
in which $h$ is a linear function in $t$ given by
 $$h(t)=(a^2-2da+d)t+a^2(d-1).$$
Since  $a<\frac{d}{2d-1}<1<d$, we have $t-a>0$ for $t\in [1,d]$, and
\begin{align*}
  h(1)=d(a-1)^2>0, \quad h(d)=(d-a)\big[d-(2d-1)a\big]>0.
\end{align*}
Therefore, $\|D-aI_K\|=d-a$, $D_a$ is positive definite, and
 $$h(t)\ge \min\big\{h(1),h(d)\big\}> 0,\quad  t\in[1,d].$$
It follows that
 \begin{align*}
   D_a-\frac{\ell(D)}{\ell(d)}= & f(D)=g(D)(dI_K-D)=A(dI_K-D)
 \end{align*}
such that $A$ is positive definite.
\end{proof}

In the following theorem, we establish criteria for the closedness of the numerical range $W(T_{{}_Q})$.

\begin{theorem}\label{thm:1stcloseT} Let  $Q\in \mathbb{B}(H)$ be a non-projection idempotent  admitted to the decomposition \eqref{equ:the unitarily of Q}, where $\widetilde{Q}$ is given by
   \eqref{eq:new decomposition of Q2 hat2} such that the operator  $D$ in \eqref{eq:new decomposition of Q2 hat2} satisfies  \eqref{equ:2propofD}.
 Let $T_{{}_Q}$ be the operator defined in \eqref{equ:nonquadopt}.   Then the following
statements are equivalent:
  \begin{enumerate}
\item[{\rm (i)}] $W(T_{{}_Q})$ is closed in $\mathbb{C}$;
\item[{\rm (ii)}] $\partial\Big(\overline{W(T_{{}_Q})}\Big)\subseteq W(T_{{}_Q})$;
\item[{\rm (iii)}] There exists $z\in \partial\Big(\overline{W(T_{{}_Q})}\Big)$ such that $z\in W(T_{{}_Q})$;
\item[{\rm (iv)}] $\ell(\|D\|)\in\sigma_p\big(\ell(D)\big)$, where $\ell(\|D\|)$ is defined in \eqref{equ:newtemp-003};
\item[{\rm (v)}] $\|D\|\in \sigma_p(D)$.
\end{enumerate}
\end{theorem}
\begin{proof}
Regarding  the equivalence of statements (i) through (v), the implications (i)$\Longrightarrow$(ii)$\Longrightarrow$(iii) are immediate, and (iv)$\Longleftrightarrow$(v) follows from Lem\-ma~\ref{lem:Du=du001}.

(iii)$\Longrightarrow$(v).
Let   $T_{{}_Q}$, $\widetilde{T_{{}_Q}}$, $D$ and $d$ be as defined  in the proof of Theorem~\ref{thm:the num of Q2}. By Theorems~\ref{thm:biggerintT} and \ref{thm:the num of Q2},
$W(\widetilde{T_{{}_Q}})=W(T_{{}_Q})$ and $\overline{W(T_{{}_Q})}$  is given  by \eqref{equ:useful003} with parameters $a,b$ and $c$ defined in \eqref{equ:abcinE}.
Let $\alpha_0$ be the angle  in  $(0,\frac{\pi}{2})$ such that
\begin{equation}\label{equ:the angle}
 \cos\alpha_0=\frac{c}{a},\quad \sin\alpha_0=\frac{b}{a}.
\end{equation}

Now suppose that $z=x+iy\in \partial\Big(\overline{W(T_{{}_Q})}\Big) \cap  W(T_{{}_Q}) $ for some $x,y\in \mathbb{R}$.
Then
\begin{equation}\label{equ:xyinbundary}
  \frac{(x-c)^2}{a^2}+\frac{y^2}{b^2}=1,
\end{equation}
and there exist $u,v\in H_5 $ with $\|u\|^2+\|v\|^2=1$ such that
$$z
=\left\langle \widetilde{T_{{}_Q}}\binom{u}{v},\binom{u}{v}\right\rangle.$$
By \eqref{equ:innertilTTT}, it follows that
\begin{align}\label{equ:formula4xy+9}
   x=\big\langle (D+I_{H_5})u,u\big\rangle-\text{Re}\Big(\big\langle \ell(D)v,u\big\rangle\Big), \quad
   y=-\text{Im}\Big(\big\langle \ell(D)v,u\big\rangle\Big).
\end{align}
If $v=0$, then from the above equations we obtain $y=0$ and $0\le x\le d+1$, which leads by $c=\frac{d+1}{2}$ to $|x-c|\le c<a$. This contradicts \eqref{equ:xyinbundary}. Therefore, $v\ne 0$. Similarly, $u\ne 0$.
Consequently, there exists  $\beta\in  (0,\frac{\pi}{2})$ such that
$$\|u\|=\cos\beta,\quad  \|v\|=\sin\beta.$$

Let $D_{-1}$ be defined by \eqref{equ:Da and SD} with  parameter $a=-1$, and let
$S'=\frac{\ell(D)}{\ell(d)}$.
Then by Lemma~\ref{lem:Da and SD}, we have
$$ I_{H_5}\ge D_{-1}\ge S'\ge 0.$$
Thus, there  exists $t\in [0,1]$ such that  $\langle D_{-1}u,u\big\rangle=t\|u\|^2$.

To express  $\big\langle \ell(D)v,u\big\rangle$, we proceed as follows:
\begin{align*}
\big|\big\langle \ell(D)v,u\big\rangle\big|
=&\ell(d)\big|\big\langle S'v,u\big\rangle\big|
=\ell(d)\big|\big\langle\big(S'\big)^{\frac12}v,\big(S'\big)^{\frac12}u\big\rangle\big|\\
\le& \ell(d)\big\|\big(S'\big)^{\frac12}v\big\|\cdot\big\|\big(S'\big)^{\frac12}u\big\|\\
\le&\ell(d)\|v\|\cdot\big\|\big(S'\big)^{\frac12}u\big\|=\ell(d)\|v\|\sqrt{\big\langle S'u,u\big\rangle}\\
\le&\ell(d)\|v\|\sqrt{\langle D_{-1}u,u\big\rangle}=\ell(d)\sqrt{t}\|v\|\cdot\|u\|\\
=&\frac{\ell(d)}{2}\sqrt{t}\sin(2\beta)=b\sqrt{t}\sin(2\beta).
\end{align*}
Therefore, there exists   $r\in[0,1]$  such that
$$\big|\big\langle \ell(D)v,u\big\rangle\big|=rb\sqrt{t}\sin(2\beta),$$
which implies
\begin{align*}
  &\text{Re}\Big(\big\langle \ell(D)v,u\big\rangle\Big)  = b\sqrt{t}r\sin(2\beta) \cos\theta, \\
  &\text{Im}\Big(\big\langle \ell(D)v,u\big\rangle\Big) =b\sqrt{t}r\sin(2\beta) \sin\theta
\end{align*}
for some $\theta\in [0,2\pi)$.  From \eqref{equ:formula4xy+9}, it follows that
\begin{align*}
  \frac{y^2}{b^2}&=tr^2\sin^{2}(2\beta) \sin^2\theta.
\end{align*}

Furthermore, using \eqref{equ:formula4xy+9}, \eqref{equ:abcinE} and \eqref{equ:the angle}, we obtain
\begin{align*}
x-c=&(d+1)\langle D_{-1}u,u\rangle-b\sqrt{t}r\sin(2\beta) \cos\theta-c\\
=&2a\cos\alpha_0\langle D_{-1}u,u\rangle-a\sin\alpha_0\sqrt{t}r\sin(2\beta) \cos\theta-a\cos\alpha_0\\
=&2at\cos\alpha_0\cos^2\beta-a\sqrt{t}r\sin\alpha_0\sin(2\beta) \cos\theta-a\cos\alpha_0.
\end{align*}
By the Cauchy--Schwarz inequality,
\begin{align*}
  \big(\frac{x-c}{a}\big)^2
=&\Big[(2t\cos^2\beta-1)\cos\alpha_0-\sqrt{t}r\sin(2\beta) \cos\theta\sin\alpha_0
\Big]^2\\
\le&(2t\cos^2\beta-1)^2+tr^2\sin^{2}(2\beta) \cos^2\theta.
\end{align*}
Therefore, if we let  $\cos^2\beta=s$, then $s\in[0,1]$ and by \eqref{equ:xyinbundary},
\begin{align*}
  1\le&(2ts-1)^2+tr^2\sin^{2}(2\beta) \cos^2\theta+tr^2\sin^{2}(2\beta) \sin^2\theta\\
   =&(2ts-1)^2+tr^2\sin^{2}(2\beta)
   \le  (2ts-1)^2+ t\sin^{2}(2\beta)\\
   =&(2ts-1)^2+4ts(1-s)=4s^2t^2-4s^2t+1\\
   \le &4s^2t^2-4s^2t^2+1=1.
\end{align*}
Hence, all inequalities, including those used in deriving
$$\big|\big\langle \ell(D)v,u\big\rangle\big|\le b\sqrt{t}\sin(2\beta),$$
are in fact equalities.  In particular,
$\langle S'u,u\rangle=\langle D_{-1}u,u\rangle$.
By Lemma~\ref{lem:Da and SD}, this implies
$$\langle A(dI_{H_5}-D)u,u\rangle=0,$$
where $A\in \mathbb{B}(H_5)$ is a positive definite operator.   It follows  that
 $Du=du$.  Therefore,  $d$ is an eigenvalue of $D$.

   (v)$\Longrightarrow$(i).  By the proof of Theorem~\ref{thm:the num of Q2}, we see that $\overline{W(T_{{}_Q})}=W(f_d)$. Since $d\in \sigma_p(D)$,  Lemma~\ref{lemma:point-spectrum-numerical-range} implies that
  $W(f_d)\subseteq W(\widetilde{T_{{}_Q}})$. Therefore, by Theorem~\ref{thm:biggerintT} we conclude
  $$\overline{W(T_{{}_Q})}\subseteq  W(\widetilde{T_{{}_Q}})= W(T_{{}_Q}),$$
   which shows that $W(T_{{}_Q})$ is closed in $\mathbb{C}$.\end{proof}

Let $Q\in \mathbb{B}(H)$ be a non-projection idempotent, and let $T_{{}_Q}$ be  the operator defined  in \eqref{equ:nonquadopt}.
 Combining Theorems~\ref{thm:the num of Q2},  \ref{thm:biggerintT} and \ref{thm:1stcloseT}, we conclude that
 $W(T_{{}_Q})$ is either a closed or an open elliptical disk. In the following corollary, we show that
for any  $r>1$, the numerical range $W(T_{{}_{Q_r}})$ is open in $\mathbb{C}$, where $Q_r$ denotes the  idempotent defined in \eqref{equ:Qr2},  and $T_{{}_{Q_r}}$ is defined in \eqref{equ:nonquadopt}  by replacing $Q$ with $Q_r$.

\begin{corollary}\label{cor:num of Qr}
For every $r> 1$, the numerical range $W(T_{{}_{Q_r}})$ is given by
\begin{equation}\label{equ:open ellipse-r}
 W(T_{{}_{Q_r}})=\Big\{x+iy:x,y\in\mathbb{R}, \frac{(x-\frac{r+3}{4})^2}{{(r^2+3r+4)}/8}+\frac{y^2}{(r^2-1)/16}<1\Big\}.
\end{equation}
\end{corollary}
\begin{proof}Let $d=\frac12(r+1)$. By Remark~\ref{rem:4identification}, Theorems~\ref{thm:the num of Q2},  \ref{thm:biggerintT} and \ref{thm:1stcloseT}, it suffices to verify that
$d$ is not an eigenvalue of the operator $D_0\in\mathbb{B}(H_0)$, where $H_0=L^2[1,d]$ and
\begin{equation}\label{equ:MD0}
  \big(D_0x\big)(t)=tx(t),\quad  x\in H_0, t\in [1,d].
\end{equation}

Suppose $u\in H_0$ satisfies $D_0u=du$. Then
$$0=\big\langle (d I_{H_0}-D_0)u,u\big\rangle=\int_{1}^{d}(d-t)|u(t)|^2dt,$$
which implies $u(t)=0$ a.e. on $[1,d]$.  Hence, $u=0$ in $H_0$, and therefore $d$ is not an eigenvalue of $D_0$.
\end{proof}

\begin{definition}For each  $r>1$, let $Q_r$ be the  idempotent defined by \eqref{equ:Qr2}. Define an operator $S_r\in C^*\{Q_r\}$ by
\begin{equation}\label{equ:Sr}
  S_r=2(1-d)\oplus 0\oplus f(D_0),
\end{equation}
where $d=\frac12(r+1)$, $D_0$ and $I_0$ are as  defined in \eqref{equ:D0I0}, and
\begin{equation*}\label{equ:W(f) of Tr}
  f(D_0)=\left(
             \begin{array}{cc}
               D_0-(2d-1)I_0 & (D_0-I_0)^\frac12(dI_0-D_0)^\frac12 \\
               0 & (2d-1)(D_0-I_0) \\
             \end{array}
           \right).
\end{equation*}
\end{definition}

\begin{theorem}\label{thm:Srnoted}
  For every $r>1$, define $S_r\in C^*\{Q_r\}$ as in \eqref{equ:Sr}. Then
  $\overline{W(S_r)}$ is not an elliptical disk.
\end{theorem}
\begin{proof} For each $t\in[1,d]$, define $f_t\in M_2(\mathbb{C})$ by
\begin{equation*}\label{equ:W(ft) of Tr}
  f_t=\left(
             \begin{array}{cc}
               t-(2d-1) & \sqrt{(t-1)(d-t)} \\
               0 & (2d-1)(t-1) \\
             \end{array}
           \right).
\end{equation*}
By  Theorem~\ref{thm:the numerical of matrix02} and Lemma~\ref{lem:the numeriacal of matrix01}, we have
\begin{equation}\label{equ:WSr02}
  \overline{W(S_r)}=\overline{\mbox{conv}\big(\cup_{t\in [1,d]}W(f_{t})\big)}
\end{equation}
such that for each $t\in [1,d]$,  $W(f_t)$ is an elliptical disk (including the  degenerate case when $t=1$ or $t=d$) with  its foci
$\lambda_1(t)$ and $\lambda_2(t)$, where $\lambda_1(t)$ and $\lambda_2(t)$ are two real numbers given by
$$\lambda_1(t)=t-(2d-1),\quad \lambda_2(t)=(2d-1)(t-1).$$
Since all these foci lie on  the $x$-axis, it follows that
$\overline{W(S_r)}$ is symmetric about the real axis; that is
$z\in  \overline{W(S_r)}$ if and only if $\bar{z}\in  \overline{W(S_r)}$ for any $z\in\mathbb{C}$.

Assume for contradiction that  $\overline{W(S_r)}$ is an elliptical disk (possibly  degenerate). Then by the above symmetry, its foci must lie on the $x$-axis. Therefore, by \eqref{equ:sf of ellipse}, there exist $x_0\in\mathbb{R}$ and $a,b\in [0,+\infty)$ such that
\begin{equation}\label{equ:hwtr-001}h_{W(S_r)}(\alpha)=x_0\cos\alpha+\sqrt{a^2\cos^2\alpha+b^2\sin^2\alpha}, \quad  \alpha\in [0,2\pi).
\end{equation}
Combining \eqref{equ:WSr02}, \eqref{equ:supportequal05} and \eqref{equ:supportequal04} yields
\begin{equation}\label{equ:same h value01}h_{W(S_r)}(\alpha)=\sup_{t\in [1,d]}h_{W(f_t)}(\alpha),\quad \alpha\in [0,2\pi).
\end{equation}

We now  compute $W(f_t)$ for each $t\in [1,d]$. By Lemma~\ref{lem:the numeriacal of matrix01} and \eqref{equ:sf of ellipse},
$$h_{W(f_t)}(\alpha)=x_t\cos\alpha+\sqrt{a_t^2 \cos^2\alpha+b_t^2\sin^2\alpha},\quad \alpha\in [0,2\pi).
$$
where
\begin{align*}
  x_t= & \frac{1}{2}\big[\lambda_1(t)+\lambda_2(t)\big]=dt-2d+1,\quad c_t=\frac{1}{2}\big|\lambda_1(t)-\lambda_2(t)\big|=(d-1)t, \\
  b_t =& \frac{\sqrt{(t-1)(d-t)}}{2},\quad a_t=\sqrt{b_t^2+c_t^2}=\frac{\sqrt{(t-1)(d-t)+\big[2(d-1)t\big]^2}}{2}.
  \end{align*}
Therefore,
\begin{align}h_{W(f_t)}(\alpha)=&(dt-2d+1)\cos\alpha\nonumber\\
\label{equ:supWBt}&+\frac12\sqrt{(t-1)(d-t)+\big[2(d-1)t\cos\alpha\big]^2}.\end{align}

Next, we consider the special cases  $t=d$ and $t=1$. Substituting  $t=d $ and $t=1$ into \eqref{equ:supWBt} gives
\begin{align}
\label{equ:hwbd}&h_{W(f_d)}(\alpha)=\begin{cases}
                      (2d-1)(d-1)\cos\alpha, & \mbox{if $\cos\alpha\ge 0$}, \\
                      (1-d)\cos\alpha, & \mbox{if $\cos\alpha< 0$},
                    \end{cases}\\
\label{equ:hwb1}&h_{W(f_1)}(\alpha)=\begin{cases}
                      0, & \mbox{if $\cos\alpha\ge 0$}, \\
                      2(1-d)\cos\alpha, & \mbox{if $\cos\alpha< 0$}.
                    \end{cases}\end{align}
Hence
$$\max\big\{h_{W(f_d)}(\alpha),h_{W(f_1)}(\alpha)\big\}=\begin{cases}
                                                          h_{W(f_d)}(\alpha), & \mbox{if $\cos\alpha\ge 0$}, \\
                                                          h_{W(f_1)}(\alpha), & \mbox{if $\cos\alpha< 0$}.
                                                        \end{cases}$$
We now  determine those   $\alpha\in [0,2\pi)$ such that
\begin{equation*}\label{equ:needed condition000}\cos\alpha\ge0 \quad \text{and}\quad h_{W(f_t)}(\alpha)\le h_{W(f_d)}(\alpha),\quad  t\in [1,d].\end{equation*}
Suppose $\cos \alpha\ge \frac{\sqrt{2}}{4\sqrt{d(2d-1)}}$. Then for any $t\in[1,d]$,
$$\cos^2\alpha\ge \frac{1}{4(2d-1)}\cdot\frac{d-1}{2d^2-d-d}\ge \frac{1}{4(2d-1)}\cdot\frac{t-1}{2d^2-d-t}.$$
Thus, for any $t\in[1,d)$,
\begin{align*}
   \cos^2\alpha\ge& \frac{(t-1)(d-t)}{4(2d-1)(2d^2-d-t)(d-t)} \\
   =&\frac{(t-1)(d-t)}{4(2d^2-d-dt)^2-4(d-1)^2t^2}.
\end{align*}
It follows that for all $t\in[1,d)$,
$$\sqrt{(t-1)(d-t)+\big[2(d-1)t\cos\alpha\big]^2}\le 2(2d^2-d-dt)\cos\alpha.$$
Note that equality holds when $t=d$. Therefore, from \eqref{equ:supWBt} and \eqref{equ:hwbd}, we conclude that
\begin{equation*}\label{equ:needed condition000} \cos\alpha\ge \frac{\sqrt{2}}{4\sqrt{d(2d-1)}}\Longrightarrow      h_{W(f_t)}(\alpha)\le h_{W(f_d)}(\alpha),\quad  t\in [1,d].\end{equation*}

Similarly,  suppose $\cos \alpha\le -\frac{\sqrt{2}}{4}$. Then for any $ t\in[1,d]$,
$$\cos^2\alpha\ge \frac18= \frac{d-1}{4\big[(2d-1)-1\big]}\ge\frac{d-t}{4\big[(2d-1)t-1\big]}.$$
Hence,  for any $ t\in[1,d]$,
\begin{align*}0\le (d-t)(t-1)\le 4\big[(2d-1)t-1\big](t-1)\cos^2\alpha,\end{align*}
which implies
\begin{align*}
0\le (d-t)(t-1)+\big[2(d-1)t\cos\alpha\big]^2\le \big[2(dt-1)\big]^2\cos^2\alpha.
\end{align*}
From \eqref{equ:supWBt} and \eqref{equ:hwb1}, it follows that
\begin{equation*}\label{equ:needed condition005}\cos \alpha\le -\frac{\sqrt{2}}{4}\Longrightarrow h_{W(f_t)}(\alpha)\le h_{W(f_1)}(\alpha),\quad  t\in [1,d].\end{equation*}

We now derive a contradiction. From the above two inequalities corresponding to the constraints on $\cos\alpha$, together with \eqref{equ:same h value01}, we conclude that
\begin{equation}\label{equ:hwtr-000}h_{W(S_r)}(\alpha)=\left\{
                                                       \begin{array}{ll}
                                                         h_{W(f_d)}(\alpha), & \hbox{if  $\alpha\in [0,\alpha_1]\cup [2\pi-\alpha_1,2\pi)$}, \\
                                                          h_{W(f_1)}(\alpha), & \hbox{if $\alpha\in [\pi-\alpha_2,\pi+\alpha_2]$,}
                                                       \end{array}
                                                     \right.
\end{equation}
where
$$\alpha_1=\mbox{arc}\cos\frac{\sqrt{2}}{4\sqrt{d(2d-1)}},\quad \alpha_2=\mbox{arc}\cos\frac{\sqrt{2}}{4}.$$
Substituting  $\alpha=0$ and $\alpha=\alpha_1$ into \eqref{equ:hwtr-001}, \eqref{equ:hwtr-000} and  \eqref{equ:hwbd}  gives
\begin{align*}&(2d-1)(d-1)=x_0+a\Longrightarrow (2d-1)(d-1)\cos\alpha_1=(x_0+a)\cos\alpha_1,\\
& (2d-1)(d-1)\cos\alpha_1=x_0\cos\alpha_1+\sqrt{a^2\cos^2\alpha_1+b^2\sin^2\alpha_1}.
\end{align*}
Therefore, $b=0$.
Then \eqref{equ:hwtr-001} implies  $h_{W(S_r)}(\pi/2)=0$. Hence,
by \eqref{equ:same h value01} and \eqref{equ:supWBt} we have
\begin{align*}0=h_{W(S_r)}(\pi/2)=\sup_{t\in [1,d]}\frac12\sqrt{(t-1)(d-t)}=\frac{d-1}{4},
\end{align*}
which contradicts  $d>1$.
\end{proof}

\begin{theorem}\label{thm:ncno}
For every $r>1$, let $S_r\in C^*\{Q_r\}$ be defined by \eqref{equ:Sr}. Then $W(S_r) $ is neither closed nor open in $\mathbb{C}$.
 \end{theorem}
\begin{proof}
We follow the notation as in the proof of  Theorem~\ref{thm:Srnoted}, and furthermore define
$$s=(2d-1)(d-1).$$
 By Theorem~\ref{thm:Srnoted}, we have
$$[x_d-a_d,x_d+a_d]=W(f_d)\subseteq \overline{W(S_r)},$$
where $$x_d=(d-1)^2, \quad a_d=d(d-1).$$
It follows that $$s=x_d+a_d\in  \overline{W(S_r)}.$$

Now, we  show that $s\notin W(S_r)$. Two observations are in order. First, from \eqref{equ:Sr} we have
\begin{equation}\label{equ:WconvWSr}
W(S_r)=\mbox{conv}\big[[2(1-d),0] \cup W\big(f(D_0)\big)\big].\end{equation}
Second, combining Definition~\ref{defn of sf}, \eqref{equ:hwtr-000}, and \eqref{equ:hwbd}, we obtain
\begin{align*}
\sup\big\{\mbox{Re}(z):z\in W(S_r)\big\}=h_{W(S_r)}(0)=s,
\end{align*}
which implies
\begin{equation}\label{x1x2+10}
x_1<s \quad \text{and} \quad \text{Re}(x_2)\le s,\quad \forall x_1\in [2(1-d),0], x_2\in W\big(f(D_0)\big).\end{equation}

Now suppose $s\in W(S_r)$. Then by \eqref{equ:WconvWSr}, there exist
$$x_1\in [2(1-d),0],\quad x_2\in W\big(f(D_0)\big),\quad t_1,t_2\in [0,1],\quad t_1+t_2=1$$
 such that $s=t_1x_1+t_2x_2$, and hence
 \begin{align*}
   s= & \text{Re}(s)=t_1x_1+t_2\text{Re}(x_2).
 \end{align*}
 If $t_1>0$, then from \eqref{x1x2+10} we deduce
 $$s<t_1s+t_2\text{Re}(x_2)\le t_1s+t_2s=s,$$ which is a contradiction.
Therefore,  $s=x_2\in W\big(f(D_0)\big)$.
Thus, there exist  $u,v\in H_{0}=L^2[1,d] $ such that
$$\|u\|^2+\|v\|^2=1\quad \text{and}\quad s=  \Big\langle f(D_0)(u,v)^T, (u,v)^T\Big\rangle.$$

Let $I_{H_0}$ be denoted simply by $I$, and write
$$f(D_0)=
\left(
  \begin{array}{cc}
     D_0-(2d-1)I & f_{12} \\
    0 & (2d-1)(D_0-I) \\
  \end{array}
\right),$$
where  $D_0\in\mathbb{B}(H_0)$ is defined by \eqref{equ:MD0},
and
$$f_{12}=(D_0-I)^{\frac12}(dI-D_0)^{\frac12}.$$
Hence, $$s=\langle D_0u,u\rangle-(2d-1)\langle u,u\rangle+\big\langle f_{12}v,u\big\rangle+\langle (2d-1)(D_0-I)v,v\rangle.$$
Since $\|u\|^2+\|v\|^2=1$, the above equation implies
\begin{equation}\label{equ:MMdplus}
(2d-1)d=\langle D_0u,u\rangle+\big\langle f_{12}v,u\big\rangle+\langle (2d-1)D_0v,v\rangle.
\end{equation}
Thus
$$(2d-1)d\le\langle D_0u,u\rangle+\big|\big\langle f_{12}v,u\big\rangle\big|+\langle (2d-1)D_0v,v\rangle.$$

Next, we derive an upper bound for $\big|\big\langle f_{12}v,u\big\rangle\big|$.
By the definition of $f_{12}$,  it is clear that
\begin{align*}
  \big|\big\langle f_{12}v,u\big\rangle\big| =& \big|\big\langle (dI-D_0)^{\frac12}v,(D_0-I)^{\frac12}u\big\rangle\big| \\
  \le&\big\|(D_0-I)^{\frac12}u\big\|\cdot
\big\|(dI-D_0)^{\frac12}v\big\|\\
\le &\frac12\Big[\big\|(D_0-I)^{\frac12}u\big\|^2+
\big\|(dI-D_0)^{\frac12}v\big\|^2\Big]\\
=&\frac12\big[\langle (D_0-I)u,u\rangle+\langle (dI-D_0)v,v\rangle\big].
\end{align*}
Substituting this bound yields
\begin{equation}\label{equ:ABuv}
(2d-1)d\le \langle Au,u\rangle+\langle Bv,v\rangle,
\end{equation}
where
\begin{align*}
   A=& D_0+\frac12(D_0-I)=\frac32D_0-\frac12I\le \frac{(3d-1)I}{2} < (2d-1)dI,   \\
   B=&(2d-1)D_0+\frac12(dI-D_0)=\big(2d-\frac32\big)D_0+\frac{1}{2}dI\le (2d-1)dI.
\end{align*}
Hence, $\langle Bv,v\rangle\le (2d-1)d\|v\|^2$, and if $u\ne 0$, then $\langle Au,u\rangle<(2d-1)d\|u\|^2$. From \eqref{equ:ABuv}, it follows that
$$(2d-1)d<(2d-1)d\|u\|^2+ (2d-1)d\|v\|^2=(2d-1)d,$$
which is a contradiction. Therefore, $u=0$ and thus $\|v\|=1$. Let
$$C=(2d-1)dI-(2d-1)D_0. $$ Then from \eqref{equ:MMdplus} we have
$\langle Cv,v\rangle= 0$.
Since $C$ is positive, this implies $Cv=0$, and thus
 $$dv-D_0v=(2d-1)^{-1}Cv=0.$$
Therefore,  $d$ is an eigenvalue of $D_0$, which is impossible as shown in  the proof of Corollary~\ref{cor:num of Qr}. This shows that  $W(S_r)$ is not closed in $\mathbb{C}$.

Finally, we prove that $W(S_r)$ is not open  in $\mathbb{C}$.
To this end, rewrite \eqref{equ:Sr} as
  $$S_r=s'\oplus 0\oplus f(D_0),$$
where $s'=2(1-d)$. From this expression, it is obvious that $s'\in W(S_r)$. Moreover, from \eqref{equ:hwtr-000} and \eqref{equ:hwb1}, we obtain
$$h_{W(S_r)}(\pi)=h_{W(f_1)}(\pi)=2(d-1)=-s'.$$
We claim that $s'\notin\text{int}\big(W(S_r)\big)$. To verify this, consider an arbitrary open disk $O_\delta=\{z\in\mathbb{C}:|z-s'|<\delta\}$ with $\delta>0$. Let $p=s'-\frac{\delta}{2}$. Then $p\in O_{\delta}$ and
$$\text{Re}\big(pe^{-i\pi}\big)=-s'+\frac{\delta}{2}>h_{W(S_r)}(\pi),$$
which implies $p\notin \overline{W(S_r)}$.
Hence, $s'\notin \text{int}\big(W(S_r)\big)$ and therefore, $W(S_r)$ is not open in $\mathbb{C}$.
\end{proof}

\end{document}